\pdfoutput=1

\newif\ifsiam
\newif\ifsubmit
\newif\ifhighlight

\siamtrue
\submittrue

\ifsiam
\documentclass[hidelinks]{siamart1116}
\else
\documentclass[a4paper]{article}
\usepackage{amsthm}
\usepackage{hyperref}
\usepackage{cleveref}
\usepackage{doi}
\usepackage[square,numbers]{natbib}
\hypersetup{pdfborder = 0 0 0, colorlinks=false,
 linkcolor=black,citecolor=black, filecolor=black, urlcolor=black, }
\newcommand{\email}{}
\fi

\usepackage{geometry}
\usepackage{amsmath,amssymb}
\usepackage{bm}
\usepackage{commath}
\usepackage{algorithm,algpseudocode}
\usepackage[title]{appendix}
\newtheorem{remark}{Remark}
\usepackage{float,placeins}

\ifhighlight
\usepackage{color,soul}
\definecolor{lightyellow}{rgb}{1.0, 1.0, 0.85}
\sethlcolor{lightyellow}
\usepackage{framed}
\newenvironment{hlbox}{%
  \MakeFramed{\advance\hsize-\width \FrameRestore}}
{\endMakeFramed}
\else
\newenvironment{hlbox}{}{}
\newcommand{\hl}[1]{#1}
\fi

\newlength\fwidth
\ifsubmit
\usepackage{tikzexternal,tikzscale}
\graphicspath{{ext}}
\else
\usepackage{tikz,tikzscale,pgfplots}
\usetikzlibrary{plotmarks}
\usetikzlibrary{external}
\pgfplotsset{compat=newest} 
\pgfplotsset{plot coordinates/math parser=false} 
\fi

\AtEndPreamble{
  \LetLtxMacro{\oldincludegraphics}{\includegraphics}  \newcommand{\myincludegraphics}[2][]{\tikzsetnextfilename{#2}\oldincludegraphics[#1]{#2}}
  \LetLtxMacro{\includegraphics}{\myincludegraphics}
}

\algdef{SE}[DOWHILE]{Do}{doWhile}{\algorithmicdo}[1]{\algorithmicwhile\ #1}%

\newcommand{\nhat}{{n}}
\newcommand{\ordo}{\mathcal{O}}
\newcommand{\Res}{\operatorname{Res}}
\newcommand{\opint}{\operatorname{I}}
\newcommand{\opquadsymb}{\operatorname{Q}}
\newcommand{\opquad}{\opquadsymb_{n}}
\newcommand{\oprem}{\operatorname{R}_{n}}

\renewcommand{\Re}{\operatorname{Re}}
\renewcommand{\Im}{\operatorname{Im}}

\newcommand{\pullb}[1]{#1}
\newcommand{\param}{\gamma}
\newcommand{\den}{\sigma}
\newcommand{\coeff}{a}
\newcommand{\upsamp}{\kappa}
\newcommand{\estq}{E_C}
\newcommand{\tol}{\epsilon}

\newcommand{\hankfk}{H^{(1)}}
\newcommand{\conj}[1]{\bar{#1}}
\newcommand{\bdry}{{\partial\Omega}}
\newcommand{\poly}{\operatorname{P}}
\newcommand{\wnum}{{k}} 
\newcommand{\work}{{W}}

\newcommand{\addAS}{s} %
\newcommand{\addAD}{d} %

\headers{Adaptive quadrature by expansion}{L. af Klinteberg and A.-K. Tornberg}
\title{Adaptive quadrature by expansion for layer potential evaluation
  in two dimensions \thanks{This work has been supported by the
    G\"oran Gustafsson Foundation for Research in Natural Sciences and
    Medicine, and by the Swedish Research Council under Grant
    No. 2015-04998.}}

\author{Ludvig af Klinteberg\thanks{Department of Mathematics /
    Swedish e-Science Research Centre (SeRC), KTH Royal Institute of
    Technology, 100 44 Stockholm, Sweden (\email{ludvigak@kth.se},
    \email{akto@kth.se})} \and Anna-Karin Tornberg\footnotemark[2]}

\begin{document}
\maketitle

\begin{abstract}
  When solving partial differential equations using boundary integral
  equation methods, accurate evaluation of singular and nearly
  singular integrals in layer potentials is crucial. A recent scheme
  for this is quadrature by expansion (QBX), which solves the problem
  by locally approximating the potential using a local expansion
  centered at some distance from the source boundary.
  In this paper we introduce an extension of the QBX scheme in 2D
  denoted AQBX -- adaptive quadrature by expansion --- which combines
  QBX with an algorithm for automated selection of parameters, based
  on a target error tolerance. A key component in this algorithm is
  the ability to accurately estimate the numerical errors in the
  coefficients of the expansion. Combining previous results for flat
  panels with a procedure for taking the panel shape into account, we
  derive such error estimates for arbitrarily shaped boundaries in 2D that
  are discretized using panel-based Gauss-Legendre quadrature.
  \hl{ Applying our scheme to numerical solutions of Dirichlet
    problems for the Laplace and Helmholtz equations, and also for
    solving these equations, we find that the scheme is able to
    satisfy a given target tolerance to within an order of magnitude,
    making it useful for practical applications.  }
  This represents a significant simplification over the original QBX
  algorithm, in which choosing a good set of parameters can be hard.
\end{abstract}



\section{Introduction}

Integral equation methods \hl{are} a class of numerical methods that
are based on the reformulation of an elliptic partial differential
equation (PDE) as a boundary integral equation. When applicable, this
solution approach has several attractive features. Among these are:
high-order discretization methods, well-conditioned linear systems
after discretization, reduced number of unknowns compared to volume
methods, and straightforward handling of moving boundaries.

A suitable starting point for our discussion on integral equation
methods is the representation of the solution $u$ using layer
potentials, which are evaluated by integrating the PDE's fundamental
solution $G$ and a layer density $\den$ over the domain boundary
$\bdry$. We represent $u$ as a linear combination of the double layer
potential $D$ and the single layer potential $S$,
\begin{align}
  u(z) = D\den(z) + \alpha S\den(z)
  = \int_\bdry \dpd{G(z,w)}{\nhat_w} \den(w) \dif s_w
  + \alpha\int_\bdry G(z, w) \den(w) \dif s_w, \quad z \in \Omega,
  \label{eq:laypot_repr}
\end{align}
\hl{ where $\nhat_w$ denotes the unit normal pointing into
  $\Omega$. For a Dirichlet problem with boundary condition $u=f$ we
  get, considering the limit $\Omega \ni z \to \bdry$ of
  \mbox{\eqref{eq:laypot_repr}} leads to a second kind integral
  equation in $\den$, }
\begin{align}
  \left( \frac{1}{2}I + D + \alpha S \right)\den(z) = f(z), \quad z \in \bdry .
  \label{eq:generic_inteq}
\end{align}
\hl{ This is a generic form for a Dirichlet problem, and we will later
  define the forms for the interior Laplace and the exterior Helmholtz
  equations.  }%
The constant $\alpha$ is selected such that the integral equation has
a unique solution and is well-conditioned \cite{Betcke2011}. Nystr\"om
discretization of this equation using a suitable quadrature rule
generates a dense linear system that can be solved rapidly by
exploiting the fact that the off-diagonal blocks are of low
rank. Solution methods include the fast multipole method (FMM)
\cite{Greengard1997} and fast direct methods
\cite{Greengard2009,Ho2012,Martinsson2005}.

The main difficulty in the procedure outlined above lies in finding a
quadrature rule that can evaluate the layer potentials $S\den(z)$ and
$D\den(z)$ when the target point $z$ is close to or on the boundary
$\bdry$. \hl{In this case,} the integrals of the layer potentials are
nearly singular or singular, requiring specialized quadrature
methods. For two-dimensional (2D) problems\hl{, several efficient such
  methods are available;} see the summaries \cite{Hao2014} and
\cite{Helsing2015}. Some of the 2D methods can be extended to
axisymmetric surfaces in three dimensions (3D)
\cite{Helsing2014}. \hl{However, for general surfaces in 3D the
  available methods are not as mature as those in 2D.} We refer to
\cite{Rahimian2016} for a recent summary of the current state of the
art.

\subsection{QBX}
\hl{In this paper, we will} focus on the relatively recent method of
quadrature by expansion (QBX) \cite{Barnett2014,Klockner2013}. The
method provides a way of evaluating both singular and nearly singular
integrals by representing the layer potential as a local expansion,
centered at a point some distance away from the boundary. While
originally proposed for Helmholtz in 2D, it can be generalized to
other PDEs \hl{in 2D and 3D; see
  \mbox{\cite{AfKlinteberg2016qbx,Askham2017,Ricketson2016}} for its
  successful use in several applications}. A strength of QBX is that
it uses the same type of local expansions as the FMM, which allows the
two methods to be combined into a fast method for evaluating layer
potentials at arbitrary locations. This is a topic of ongoing research
\cite{Rachh2017}. As with most methods, however, QBX solves one
problem and introduces another. \hl{While the problem of evaluating
  layer potentials on or very close to the boundary is solved, it is
  replaced with} the new problem is how to efficiently compute the
local expansion of the layer potential. In particular, computing the
expansion coefficients entails evaluating a series of integrals with
increasing order of near singularity. \hl{Effectively one has} traded
\hl{a} hard problem for several easier problems. Nevertheless, QBX is
still a competitive method for these problems, especially since it has
a solid analytical groundwork \cite{Epstein2013}.

One of the difficulties of QBX is that of parameter selection. The
convergence of the local expansion is governed by the order $p$ at
which it is truncated, and by the distance $r$ between the boundary
and the expansion center. The expansion coefficients are computed
using a quadrature rule for smooth integrands, and for them to be
accurate it is necessary to upsample the boundary points by some
factor $\kappa$. These three parameters together affect the two
competing errors of QBX: the truncation error and the coefficient
error (often referred to as the quadrature error). The truncation
error increases with $r$ and decreases with $p$, while the coefficient
error increases with $p$ and decreases with $\upsamp$ and $r$ (see
\cite[Fig. 3]{AfKlinteberg2016qbx} for an example). Together $r$, $p$
and $\kappa$ form a large parameter space, and how to best set these
parameters is not clear. Instead, experimentation must be used to
determine good parameter ranges for a specific application. This is in
itself not unfeasible, but it would be preferable to reduce the number
of free parameters.

\subsection{Contribution}

In this paper we propose a scheme for adaptively setting the order $p$
and upsampling \hl{factor} $\upsamp$ at the time of computation, such
that a the error is maintained below a target tolerance. The key
ingredient for this \hl{scheme} to be successful is the ability to
accurately estimate the magnitude of the coefficient error, which is
the quadrature error in the expansion coefficients. Such estimates
were derived in \cite{AfKlinteberg2016quad} for simple geometries in
two dimensions\hl{, namely} flat Gauss-Legendre panels and the
trapezoidal rule on the unit circle. Here we extend these estimates,
by locally using a polynomial to represent the mapping between a flat
panel and a panel of general shape. This greatly increases the
accuracy of the estimates, and allows us to build an adaptive QBX
scheme based on them. We here restrict ourselves to analyzing
Gauss-Legendre panel quadrature, but the methodology could equally
well be applied to a discretization using the trapezoidal rule.

\hl{Taking into account the mapping of the parametrization when
  analyzing nearly singular quadrature errors is by itself not new;} a
discussion on nearly singular quadrature errors similar to ours can be
found in \cite{Barnett2014}. Our main contribution in this regard is
the construction of an explicit representation of the mapping, which
we then invert in order to compute an error estimate.

This paper is organized as follows: In \cref{sec:foundations-aqbx} we
introduce the general structure of our scheme, and derive the details
for the Laplace double layer potential and the Helmholtz combined
field potential. In \cref{sec:quadrature-errors} we show how to
compute the coefficient error estimates necessary for the scheme to be
useful. In \cref{sec:local-aqbx} we briefly discuss how to combine the
scheme with a fast method. In \cref{sec:numer-exper} we show a
selection of numerical results \hl{that} illustrate the performance of
our method.

\section{Foundations of AQBX}
\label{sec:foundations-aqbx}

In this section we begin by introducing the foundations of our method
using a generic notation, before we give the specific details for the
Laplace and Helmholtz equations. We start from a given layer potential
representation,
\begin{align}
  u(z) = \int_\bdry G(z, w) \den(w) \dif s_w .
\end{align}
To evaluate this using AQBX, we first need to split the fundamental
solution using a suitable addition theorem, \newcommand{\addA}{A}
\newcommand{\addB}{B}
\begin{align}
  G(z, w) = \sum_{m=0}^\infty \addA_m^r(w, z_0) \cdot \addB_m^r(z, z_0) ,
  \label{eq:add_thm_gen}
\end{align}
where $A_m^r$ and $B_m^r$ are either vectors or scalars, depending on
the fundamental solution.  For a given $r$, the above addition theorem
is valid for
\begin{align}
  |z-z_0| < r \le |w-z_0|.
\end{align}
For the specific formulas for Laplace and Helmholtz, see
\eqref{eq:cauchy_taylor} and \eqref{eq:grafadd}.  \hl{We assume that
  \mbox{\eqref{eq:add_thm_gen}} is normalized such that the following
  holds for $|z-z_0| \le r$ and integer $m\ge 0$:}
\begin{hlbox}
  \begin{align}
    \begin{split}
      \abs{\addB_m^r(z, z_0)} &\le 1 , \\
      \max_{m, z} \abs{\addB_m^r(z, z_0)} &= 1.
    \end{split}
                                            \label{eq:normalization}
  \end{align}
\end{hlbox}
Here, and throughout the paper, we let $|\cdot|$ denote the
$\ell^2$-norm for vectors, and the complex modulus for scalars.
\hl{If the second condition holds for each $m$, i.e.
  $\max_{z} \abs{\addB_m^r(z,z_0)} = 1$, then we get slightly sharper
  error estimates in \mbox{\cref{alg:coeff}}. This is however not
  necessary, and indeed only holds for the Laplace formulation used in
  this paper.} %
The addition theorem allows us to pick an expansion center $z_0$,
determine $r$ as
\begin{align}
  r = \min_{w \in \bdry} |w-z_0|,
\end{align}
and evaluate the layer potential as a local expansion in the
neighborhood of $z_0$,
\begin{align}
  u(z) = \sum_{m=0}^\infty \coeff_m \cdot \addB_m^r(z, z_0), \quad |z-z_0| \le r,
  \label{eq:local_exp}
\end{align}
where
\begin{align}
  \coeff_m = \int_\bdry \addA_m^r(w, z_0) \den(w) \dif s_w .
  \label{eq:int_coeff}
\end{align}
The fact that \eqref{eq:local_exp} holds also at the equality
$|z-z_0|=r$ was shown in \cite{Epstein2013} for the Laplace and
Helmholtz equations, but can be generalized to other kernels,
e.g. Stokes equations \cite{AfKlinteberg2016qbx}. %
\hl{ This allows us to evaluate the expansion also on $\bdry$, at the
  single point which is closest to $z_0$. In fact, in
  \mbox{\cite{Barnett2014}} it is shown that the expansion even
  converges in a disc about $z_0$ with radius $R > r$, as long as the
  density $\sigma$ is analytic inside that disc.}

In a computational scheme the local expansion is truncated at a
maximum order $p$, and the coefficients, which we now denote
$\tilde\coeff_m$, are computed using a suitable quadrature
rule. This gives us the QBX approximation of the layer potential,
\begin{align}
  u_p(z) = \sum_{m=0}^p \tilde\coeff_m \cdot \addB_m^r(z, z_0) .
  \label{eq:local_exp_p}
\end{align}
The error in this approximation can, \hl{assuming exact arithmetic},
be separated into a truncation error $e_T$ and a coefficient error
$e_C$ \cite{Epstein2013},
\begin{align}
  u(z) - u_p(z) =
  \underbrace{u(z) - \sum_{m=0}^p \coeff_m \cdot \addB_m^r(z, z_0)}_{e_T}
  + 
  \underbrace{\sum_{m=0}^p (\coeff_m - \tilde\coeff_m) \cdot \addB_m^r(z, z_0)}_{e_C} .
  \label{eq:err_split}
\end{align}

\subsection{Truncation error}
\label{sec:truncation-error}

The truncation error $e_T$ of the scheme arises because we limit the
local expansion to $p+1$ terms. In our normalized form
\eqref{eq:normalization} it can be bounded by
\begin{align}
  \abs{e_T} = \abs{\sum_{m=p+1}^\infty \coeff_m \cdot \addB_m^r(z, z_0)} 
  \le \sum_{m=p+1}^\infty \abs{\coeff_m} .
\end{align}
It was shown in \cite{Epstein2013} that the truncation error for
several kernels satisfies
\begin{align}
  e_T \le M(p, r) r^{p+1} \norm{u}_{\mathcal C^{p+1}(|z-z_0|\le r)},
  \label{eq:eT_bound}
\end{align}
\hl{for some positive $M(p, r)$. In our experience, the truncation
  error typically decays exponentially in $p$, with a rate that is
  proportional to $r$. Finding a usable a priori estimate for the
  error is hard, since it depends on both the local geometry of
  $\bdry$ and the regularity of the density $\den$ in a nontrivial
  way. However, assuming that the expansion coefficients decay
  exponentially, a good estimate for the truncation error is }
\begin{align}
  \abs{e_T} \approx \abs{ \coeff_{p+1} \cdot \addB_{p+1}^r(z, z_0) }
  \le \abs{ \coeff_{p+1} } .
\end{align}
In practice we only have the coefficients up to $\coeff_p$. Therefore,
we can define a useful (and usually conservative) a posteriori
estimate as
\begin{align}
  \abs{e_T} \approx \abs{ \coeff_p \cdot \addB_p^r(z, z_0) }
  \le \abs{ \coeff_p } .
\end{align}

\subsection{Coefficient error}

The coefficient error $e_C$ in \eqref{eq:err_split} is a result of the
numerical evaluation of the coefficient integrals \eqref{eq:int_coeff}
for $m=0,\dots,p$, which for a given boundary $\bdry$ is evaluated
using some $n$-point quadrature rule. Assuming that the density $\den$
and the boundary $\bdry$ are well-resolved by that quadrature (which is
a prerequisite for the underlying Nystr\"om discretization), the main
source of the error is the near singularity in $\addA_m^r$ when
evaluated at $z_0$. The order of this near singularity typically grows
with $m$, and to counter this the density must be upsampled (by
interpolation) to a grid which is fine enough to resolve
$\addA_m^r$. The amount of upsampling required can be determined
through $\estq(n, m)$, which is an accurate a priori estimate of the
coefficient error in term $m$ when using $n$ quadrature nodes,
\begin{align}
  \estq(n, m) \approx \abs{\coeff_m - \tilde\coeff_m} \ge 
  \abs{(\coeff_m - \tilde\coeff_m) \cdot \addB_m^r(z, z_0)} .
\end{align}
We will in section \ref{sec:quadrature-errors} show how to derive such
error estimates for the cases of the Laplace and Helmholtz equations.

\begin{remark}[Resolution error]
  Here we only discuss errors in the coefficients $\tilde\coeff_m$ due
  to near singularities in the integrals
  \mbox{\eqref{eq:int_coeff}}. However, there is also a lower bound on
  the accuracy of the coefficients, imposed by how well the underlying
  grid resolves the boundary $\bdry$ and the layer density $\den$. For
  a given panel, this error could be estimated by analyzing a modal
  expansion of the grid point coordinates and density values, such as
  the Legendre polynomial expansion used in
  \cref{sec:quadrature-errors}. \hl{This is briefly explored in
    \mbox{section~\ref{sec:source-points-close}}, though a more
    in-depth analysis is beyond the scope of the present paper.}
\end{remark}

\subsection{The AQBX scheme}

Let us now assume that we have a discretization of $\bdry$
characterized by $n$, which for a global quadrature denotes the total
number of points, and for a panel-based quadrature denotes the number
of points on each panel. Denoting that quadrature $\opquad$, we define
a combined interpolation and quadrature operator
$\opquadsymb_{\upsamp n}$, which computes the quadrature by first
upsampling the density to $\upsamp n$ points (for simplicity we assume
$\upsamp$ integer). Furthermore, we assume that the error when
computing a coefficient $a_m$ is well estimated by a function
$\estq(\upsamp n, m)$. Then the adaptive algorithm for evaluating
$u(z)$ in the neighborhood of $z_0$ to a tolerance $\epsilon$ can be
summarized using \cref{alg:coeff,alg:eval}.
\begin{algorithm}[htbp]
  \caption{Compute expansion coefficients at $z_0$ to tolerance $\tol$.}
  \label{alg:coeff}
  \begin{algorithmic}
    \Function{Expansion coefficients}{$z_0, \tol$}
    \State $\upsamp \gets 1, \: m \gets 0$
    \Repeat
    \While{$\estq(\upsamp n, m) > \tol$}
    \Comment{Check $e_C$ estimate}
    \State $\upsamp \gets \upsamp+1$
    \Comment{Increase upsampling rate}
    \EndWhile
    \State $a_m \gets \opquadsymb_{\kappa n}[\addA_m^r(\cdot, z_0) \den(\cdot)]$
    \Comment{Compute coefficient}
    \State $m \gets m+1$
    \Until{$\abs{a_m} < \tol$}
    \Comment{Break when $e_T$ estimate below tolerance}
    \State \Return \hl{$\{a_m\}_{m=0}^p$}
    \EndFunction
  \end{algorithmic}
\end{algorithm}
\begin{algorithm}[htbp]
  \caption{Evaluate $u(z)$ to tolerance $\tol$ using expansion at $z_0$.}
  \label{alg:eval}
  \begin{algorithmic}
    \Function{Evaluate expansion}{$z, z_0,$ \hl{$\{a_m\}_{m=0}^p$}, $\tol$}
    \State $u \gets 0, \: m \gets 0$
    \Repeat
    \State $\delta_m \gets \coeff_m \cdot \addB_m^r(z, z_0)$
    \Comment{Evaluate $m$th  expansion term}
    \State $u \gets u + \delta_m$
    \State $m \gets m+1$
    \Until{$\abs{\delta_m} < \tol$}
    \Comment{Break when $e_T$ estimate below tolerance}
    \State \Return{$u$}
    \EndFunction
  \end{algorithmic}
\end{algorithm}
Note that if \cref{alg:coeff} has produced $p+1$
coefficients, then \cref{alg:eval} is guaranteed to terminate
within $p+1$ iterations, since $|\delta_p| \le |a_p| < \tol$. For more
conservative termination criteria, one can use
$\max(|a_{m-1}|,|a_m|) < \tol$ and
$\max(|\delta_{m-1}|,|\delta_m|) < \tol$ to take into account any
even/odd behavior in the expansion. 

\subsection{Specific applications}
\label{sec:aqbx_examples}

We will now proceed with formulating AQBX for two different
applications: the Laplace double layer potential and the Helmholtz
combined field potential.

\subsubsection{Laplace equation}
\label{sec:laplace-equation}

We first consider the Laplace Dirichlet problem in a domain $\Omega$ bounded by a boundary $\bdry$,
\begin{align}
  \Delta u &= 0, \quad \mbox{in } \Omega, \\
  u &= f, \quad \mbox{on } \bdry.
\end{align}
The fundamental solution to this PDE is
\begin{align}
  \phi(z,w) = \frac{1}{2\pi} \log\abs{z-w} .
\end{align}
The interior Dirichlet problem can be represented using the double
layer potential,
\begin{align}
  u(z) = D\den(z) = \int_\bdry \dpd{\phi(z,w)}{\nhat_w} \den(w) \dif s_w .
  \label{eq:laplace_dbl_lyr}
\end{align}
In complex notation this can be compactly represented as
\begin{hlbox}
\begin{align}
  u(z) &= \Re v(z), \\
  v(z) &= \frac{1}{2\pi} \int_\bdry \frac{\nhat_w}{z-w} \den(w) \dif s_w . 
\end{align}
\end{hlbox}
On $\bdry$ the density $\den$ satisfies the integral equation
\begin{align}
  \left( \frac{1}{2}I + D \right) \den = f .
  \label{eq:laplace_inteq}
\end{align}
It should be noted that this particular layer potential actually has a
smooth limit on the boundary, so no special quadrature is needed for
solving the integral equation, unless different parts of the boundary
are close to each other. Special quadrature is however still needed
for evaluating the solution close to the boundary, once $\den$ has
been computed.

Starting from the Taylor expansion of the Cauchy kernel,
\begin{align}
  \frac{-1}{z-w} = \sum_{m=0}^\infty \frac{(z-z_0)^m}{(w-z_0)^{m+1}},
  \label{eq:cauchy_taylor}
\end{align}
it is straightforward to derive the AQBX formulation
\eqref{eq:add_thm_gen} for $v(z)$,
\begin{hlbox}
\begin{align}
  \addA_m^r(w, z_0) &= -\frac{r^m \nhat_w}{2\pi (w-z_0)^{m+1}}, \\
  \addB_m^r(z, z_0) &= \frac{(z-z_0)^m}{r^m} .
\end{align}
\end{hlbox}
Note that, by these definitions together with \eqref{eq:int_coeff},
the real part of the coefficient $\coeff_0$ will simply hold the value
of the double layer potential evaluated at $z_0$, i.e. 
$\Re \coeff_0 = u(z_0)$.

\subsubsection{Helmholtz equation}
\label{sec:helmholtz-equation}

We now consider the Helmholtz Dirichlet problem in an unbounded domain
$\Omega$ exterior to a boundary $\bdry$, which for a wavenumber $k$
is stated as
\begin{align}
  \Delta u + k^2 u &= 0, \quad \mbox{in } \Omega,\\
  u &= f, \quad \mbox{on } \bdry .
\end{align}
The solution must satisfy the Sommerfeld radiation condition for $r = |z|$,
\begin{align}
  \label{eq:sommerfeld}
  \lim_{r \to \infty} r^{1/2}\left( \dpd{u}{r} - i k u \right) = 0,
\end{align}
which gives a fundamental solution that is essentially the
zeroth-order Hankel function of the first kind,
\begin{align}
  \label{eq:fundsol}
  \phi_k(z, w) = \frac{i}{4} \hankfk_0(k|z-w|) .
\end{align}
It is possible to represent the solution using the combined field
integral representation\footnote{The general form of the
  representation is $D_k - i\eta S_k$. We have here set $\eta = k/2$,
  following \cite{Helsing2015}.},
\begin{align}
  \label{eq:combfield}
  u(z) = \int_\bdry \left(\dpd{\phi_k(z, w)}{n_w} - 
  \frac{ik}{2} \phi_k(z, w) \right) \den(w) \dif s_w
  = \left(D_k\den - \frac{ik}{2} S_k\den \right) (z) .
\end{align}
Here the double layer kernel is
\begin{align}
  \dpd{\phi_k(z, w)}{n_w}
  &= \frac{ik}{4} \hankfk_1(k|z-w|) \frac{(z-w) \cdot n_w}{|z-w|}  \\
  &= \frac{ik}{4} \hankfk_1(k|z-w|) |z-w| \Re \left[ \frac{n_w}{z-w} \right] .
    \label{eq:helmholtz_dbl_kern}
\end{align}
On the boundary we get the integral equation
\begin{align}
  \left( \frac{1}{2}I + D_k - \frac{ik}{2} S_k \right) \den = f .
  \label{eq:helmholtz_inteq}
\end{align}

To formulate AQBX for the combined field representation, we start from
the Graf addition theorem \cite[\S10.23(ii)]{NIST:DLMF},
\begin{align}
  \frac{i}{4} \hankfk_0(k|z-w|) 
  &= \sum_{m=-\infty}^\infty 
    \frac{i}{4} \hankfk_m(kr_w) e^{-im\theta_w}
    J_m(kr_z) e^{im\theta_z}, 
    \quad r_z < r_w .
    \label{eq:grafadd}
\end{align}
Here $(r_w,\theta_w)$ and $(r_z,\theta_z)$ are the polar coordinates
of $w-z_0$ and $z-z_0$,
\begin{align}
  \begin{aligned}
    r_w &= |w-z_0|, & r_z &= |z-z_0|, \\
    e^{-im\theta_w} &= \frac{|w-z_0|^m}{(w-z_0)^m}, \quad & e^{im\theta_z} &= \frac{(z-z_0)^m}{|z-z_0|^m}. 
  \end{aligned} \label{eq:exppot}
\end{align}
\begin{hlbox}
  We can thus form an expansion for the kernel of the combined field
  representation \eqref{eq:combfield} as 
  \begin{align}
    \dpd{\phi_k(z, w)}{n_w} - \frac{ik}{2} \phi_k(z, w) &=
    \sum_{m=-\infty}^\infty
    c_m(w, z_0)
    J_m(kr_z) e^{im\theta_z},     \label{eq:combfield_expansion} \\
    c_m(w, z_0) &=  \addAD_m(w, z_0) - \frac{ik}{2} \addAS_m(w, z_0)  .
  \end{align}
  Here $\addAS_m$ is an immediate result of \eqref{eq:grafadd},
  \begin{align}
    \addAS_m(w, z_0) &= \frac{i}{4} \hankfk_m(kr_w) e^{-im\theta_w},
                       \label{eq:addAS}
  \end{align}
  while $\addAD_m$ is obtained by differentiation of \eqref{eq:addAS}. A
  compact form for $\addAD_m$ was derived in \cite{Barnett2014},
  \begin{align}
    \addAD_m(w, z_0) &= \frac{ik}{8} \left(
                       \hankfk_{m-1}(kr_w) e^{-i(m-1)\theta_w} \conj n_w
                       -\hankfk_{m+1}(kr_w) e^{-i(m+1)\theta_w} n_w
                       \right) .
                       \label{eq:addAD}
  \end{align}  
The Bessel functions $J_m$ decay rapidly with $m$. It however easy to
experimentally verify that a normalization satisfying
\eqref{eq:normalization} is obtained by using the first term of the
power series for $J_m(kr)$ \cite[\S10.4,\S10.8]{NIST:DLMF},
\begin{align}
  J_{\pm m}(kr) = (\pm 1)^m \frac{1}{m!} \left( \frac{kr}{2} \right)^m
  + \ordo\left( \frac{1}{(m+1)!} \left( \frac{kr}{2} \right)^{m+2} \right),
  \quad m \ge 0.
  \label{eq:bessel_norm}
\end{align}
We have defined our general AQBX formulation \eqref{eq:add_thm_gen}
for indices $m\ge 0$. To fit the expansion
\eqref{eq:combfield_expansion} into this we let $\addA_m^r$ and
$\addB_m^r$ be vector-valued for $m>0$, representing the terms with
indices $\pm m$ in \eqref{eq:combfield_expansion}. The AQBX
formulation for the combined field representation is then given by
\begin{align}
  \addA_0^r(w, z_0) &= c_0(w, z_0),
                      \label{eq:addA0_helm}\\
  \addB_0^r(z, z_0) &= J_0(kr_z),
                      \label{eq:addB0_helm}.
\end{align}
and, for $m>0$,
\begin{align}
  \addA_m^r(w, z_0) &= \frac{\sqrt{2}}{m!}  \left( \frac{kr}{2} \right)^{m} 
                      \left( c_m(w, z_0), \: c_{-m}(w, z_0) \right),
                      \label{eq:addA_helm}\\
  \addB_m^r(z, z_0) &= \frac{m!}{\sqrt{2}} \left( \frac{2}{kr} \right)^{m} 
                      \left( J_m(kr_z) e^{im\theta_z}, \: J_{-m}(kr_z) e^{-im\theta_z} \right)
                      \label{eq:addB_helm},
\end{align}
such that the coefficients $a_m \in \mathbb C^2$ for $m>0$.
\end{hlbox}

We again note that the coefficient $a_0$ \eqref{eq:int_coeff} will
hold the value of the potential at $z_0$, since (trivially)
$\addAS_0 = \phi_k$ in \eqref{eq:fundsol} and (through
\eqref{eq:helmholtz_dbl_kern})
\begin{equation}
    \addAD_0(w, z_0) =
    -\frac{ik}{4} \hankfk_{1}(kr_w) \Re \left[
      e^{-i\theta_w} n_w
    \right] 
    = \dpd{\phi_k(z_0, w)}{n_w} .
  \label{eq:addADzero}
\end{equation}

\section{Coefficient errors}
\label{sec:quadrature-errors}

In this section we derive a central piece of AQBX: the coefficient
error estimate $\estq$ required in \cref{alg:coeff}. We will consider
the layer potentials of section \ref{sec:aqbx_examples} combined with
a panel-based quadrature, where the boundary curve is subdivided
into smaller segments, each of which is discretized using an $n$-point
Gauss-Legendre quadrature. This is sometimes referred to as a
composite Gauss-Legendre quadrature.

We begin by considering the error in the contribution to a coefficient
from a single panel; the total error can then be computed as the sum
of errors from all adjacent panels. %
\hl{ For this, let $\Gamma$ be an open curve (i.e. a panel)
  parametrized by an analytic function $\param(t) \in \mathbb{C}$,
  $t\in[-1,1]$, with the normal defined as
  $\nhat(t)=i\param'(t) / |\param'(t)|$, and $\gamma$ oriented such
  that $\nhat$ points into the domain $\Omega$.} %
For simplicity we denote by $\den(t)$ the pullback of $\den$ under
$\gamma$, i.e.  $\pullb\den(t) = \den(\param(t))$.

\subsection{Laplace coefficient error}
Beginning with QBX for the Laplace double layer potential
(sec. \ref{sec:laplace-equation}), the expansion coefficients are
computed as
\begin{align}
  \coeff_m 
  &= -\frac{r^m}{2\pi} \int_\Gamma \frac{\den(w) n_w}{(w-z_0)^{m+1}} \dif s_w \label{eq:dbl_lyr_coef} \\
  &= -\frac{ir^m}{2\pi} \int_{-1}^1 \frac{\pullb\den(t)}{(\param(t)-z_0)^{m+1}} \param'(t) \dif t .
    \label{eq:dbl_lyr_coef_t}
\end{align}
The coefficient errors are introduced when this integral is evaluated
using a discrete quadrature rule. In a boundary integral equation
context, we assume a panel $\Gamma$ such that $\param$ and
$\pullb\den$ are well-resolved by an $n$-point Gauss-Legendre
quadrature rule. We will here focus on the standard choice
$n=16$.

\begin{hlbox}
  For our discussion on quadrature errors, we introduce the following
  notation: Let $\opint$ denote the integral over $[-1, 1]$, and let
  $\opquad$ denote the Gauss-Legendre approximation of that integral,
\begin{align}
  \opint[f] &= \int_{-1}^1 f(x) \dif x, \\
  \opquad[f] &= \sum_{i=1}^n f(x_i) w_i,
\end{align}
The quadrature error $\oprem$ is then defined as
\begin{align}
  \oprem[f] = \opint[f] - \opquad[f],
  \label{eq:oprem}
\end{align}
and $\opint$, $\opquad$ and $\oprem$ are all linear functionals on
$C(-1,1)$.
\end{hlbox}

If the above assumptions on the resolution hold, then we can expect
the quadrature error to be small for the integrand
\eqref{eq:dbl_lyr_coef_t}, provided that $z_0$ is far away from
$\Gamma$. If on the other hand $z_0$ is close to $\Gamma$, then the
quadrature error will be dominated by the nearly singular quadrature
error that arises when the integrand is evaluated close to its
pole.

\begin{hlbox}
  To estimate nearly singular quadrature errors, we can can proceed in
  the same way as in \cite{AfKlinteberg2016quad}. The central property
  which we will use is the following: Let $f(t)$ be a function which
  is analytic on $[-1, 1]$ and everywhere inside a contour
  $\mathcal C$ enclosing $[-1,1]$, except at a finite number of poles
  $\{t_j\}_{j=0}^{N}$ enclosed by $\mathcal C$. If $f$ is integrated
  on $[-1,1]$ using the $n$-point Gauss-Legendre rule, then the
  quadrature error \eqref{eq:oprem} is given by \cite{Donaldson1972}
  \begin{align}
    \oprem[f] = \frac{1}{2\pi i}\int_{\mathcal C} f(t)k_n(t) \dif t  -
    \sum_{j=0}^{N} \Res\left[ f(t) k_n(t), t_j \right],
    \label{eq:quaderr_contour}
  \end{align}
  where $k_n$ is the characteristic remainder function,
  \begin{align}
    k_n(t) = \frac{1}{P_n(t)} \int_{-1}^1 \frac{P_n(s)}{t-s}\dif s.
  \end{align}
  While not available in closed form, it can in the limit $n\to\infty$
  be shown to satisfy \cite[appx.]{Barrett1961}
  \begin{align}
    k_n(t) = \frac{2\pi}{(t \pm \sqrt{t^2-1})^{2n+1}} \left( 1 + \ordo\left(1/n\right) \right),
    \label{eq:kn}
  \end{align}
  with the sign in the denominator given by the sign of $\Re t$. This
  result, though asymptotic, is accurate already for small $n$
  \cite{AfKlinteberg2016quad}.  Since $n \ge 1$, we can bound the
  remainder function as
  \begin{align}
    \abs{k_n(t)} \le \frac{2 \pi C(t)}{|t \pm \sqrt{t^2-1}|^{2n+1}},
  \end{align}
  for some $C(t) > 0$. If $t$ lies on the Bernstein ellipse $B_\rho$,
  defined as the ellipse with foci $\pm 1$ where the semimajor and
  semiminor axes sum to $\rho>1$, then $|t \pm \sqrt{t^2-1}|=\rho$ and
  \begin{align}
    \abs{k_n(t)} \le \frac{2 \pi C_\rho}{\rho^{2n+1}},
  \end{align}
  where $C_\rho = \max_{t\in B_\rho}C(t)$. If $|f| \le M_\rho$ on some
  Bernstein ellipse $B_\rho$, and we let the contour integral in
  \eqref{eq:quaderr_contour} follow that ellipse, then the integral
  can be bounded as
  \begin{align}
    \abs{ \frac{1}{2\pi i}\int_{B_\rho} f(t)k_n(t) \dif t }
    \le \frac{C_\rho M_\rho}{\rho^{2n+1}} \int_{B_\rho} |\dif t|
    \le \frac{4 C_\rho M_\rho}{\rho^{2n}} .
    \label{eq:contour_bound}
  \end{align}
  The last above step is motivated by the geometry of a Bernstein
  ellipse; the circumference of $B_\rho$ has its maximum as
  $\rho\to 1$, in which case the ellipse collapses onto $[-1,1]$. The
  contribution from the contour integral in \eqref{eq:quaderr_contour}
  vanishes completely if $f(t) k_n(t) t \to 0$ as $|t| \to \infty$,
  and is generally small compared to the residue for poles close to
  $[-1,1]$, as was noted in \cite{AfKlinteberg2016quad}.

  In the case where $f$ has a simple pole $t_0$, we can bound the
  error as
  \begin{align}
    \abs{\oprem[f]} \le \frac{2 \pi C(t_0)}{|t_0 \pm \sqrt{t_0-1}|^{2n+1}}
    \lim_{t\to t_0} \abs{(t-t_0)f(t_0)} + \frac{4 C_\rho M_\rho}{\rho^{2n}}.
    \label{eq:f_simple_bound}
  \end{align}
  Now consider the more general case, when $f$ has a single pole $t_0$
  of order $m+1$ (the case of several poles follows trivially).  For
  the residue we can write
  \begin{align}
    \Res\left[ f(t) k_n(t), t_0 \right] = \frac{1}{m!}
    \sum_{\ell=0}^m {\ell\choose m} k_n^{(m-\ell)}(t_0) \lim_{t\to t_0}
    \left[ \dod[\ell]{}{t} (t-t_0)^{m+1} f(t) \right],
    \label{eq:res_splitsum}
  \end{align}
  since $k_n(t_0)$ is smooth for $t \notin [-1,1]$. Defining
  \begin{align}
    \norm{f}_{C^m(t_0)} := \max_{\ell \le m} \lim_{t\to t_0} \abs{\dod[\ell]{f}{t}},
  \end{align}
  the residue can be bounded as
  \begin{align}
    \abs{\Res\left[ f(t) k_n(t), t_0 \right]} \le \frac{2^m}{m!}
    \big\|k_n\big\|_{C^m(t_0)} \big\| (t-t_0)^{m+1} f(t) \big\|_{C^m(t_0)} .
  \end{align}
  This allows us to bound the quadrature error as
  \begin{align}
    \abs{\oprem[f]} \le \frac{2^m}{m!}
    \big\|k_n\big\|_{C^m(t_0)} \big\| (t-t_0)^{m+1} f(t) \big\|_{C^m(t_0)}
    + \frac{4 C_\rho M_\rho}{\rho^{2n}}.
    \label{eq:f_general_bound}
  \end{align}  
  
  The above bound will overestimate the error by a large factor. In
  practice a good approximation of the error is achieved if we neglect
  the contribution from the contour, under the assumption that the
  pole is close to the interval, and assume that $k_n$ varies much
  more rapidly with $t$ than the other factors. Only keeping the term
  in \eqref{eq:res_splitsum} with the highest derivative in $k_n$, we
  get the error approximation
  \begin{align}
    \oprem[f] &\approx -  \frac{1}{m!} k_n^{(m)}(t_0)
                \lim_{t \to t_0} \left( (t-t_0)^{m+1}
                f(t)  \right).
                \label{lem:quad_err} 
  \end{align}
  This, with the remainder function evaluated using \eqref{eq:kn}, is
  the approximation that we will be using henceforth. The derivatives
  of $k_n$ can (at least for small $m$) be derived analytically from
  \eqref{eq:kn}, or, as shown in \cite{AfKlinteberg2016quad},
  approximated as
  \begin{align}
    k_n^{(m)}(t) \approx k_n(t) \left( \mp \frac{2n+1}{\sqrt{t^2-1}}   \right)^m .
    \label{eq:knm}
  \end{align}

  Returning to the QBX coefficients of the Laplace double layer
  potential, the integral \eqref{eq:dbl_lyr_coef} clearly has a pole of
  order $m+1$ at $z_0$. In parametrized form \eqref{eq:dbl_lyr_coef_t}
  the pole instead lies at the point $t_0 \in \mathbb C$, such that
  \begin{align}
    \param(t_0) - z_0 = 0.
    \label{eq:t0_def}
  \end{align}
  Note that evaluating $\param(t_0)$ corresponds to evaluating an
  analytic continuation of $\gamma(t)$ in some neighborhood of
  $[-1,1]$.  Denoting by $\tilde a_m$ the approximation of $a_m$ (as
  defined in \eqref{eq:dbl_lyr_coef_t}), and also assuming that there
  exists an analytic continuation of $\pullb\den$, then we have from
  \eqref{lem:quad_err} that the error can be approximated as
  \begin{align}
    \coeff_m - \tilde \coeff_m 
    = \oprem\left[ -\frac{ir^m \pullb\den \param'}{2\pi (\param-z_0)^{m+1}} \right]
    \approx
    \frac{ir^m}{2\pi m!} k_n^{(m)}(t_0) \lim_{t \to t_0} 
     \frac{(t-t_0)^{m+1} \pullb\den(t) \param'(t)}{(\param(t)-z_0)^{m+1}}  .
    \label{eq:laplace_coeff_err}
  \end{align}
  This we can simplify by noting that since $z_0=\param(t_0)$,
  \begin{align}
    \lim_{t \to t_0} \frac{(t-t_0)^{m+1}}{(\gamma(t) - z_0)^{m+1}} = \frac{1}{(\gamma'(t_0))^{m+1}} .
    \label{eq:limit_gives_deriv}
  \end{align}
  Thus a coefficient error estimate
  $\estq(n, m)\approx|\coeff_m-\tilde\coeff_m|$ is given by combining
  \eqref{eq:kn}, \eqref{eq:knm}, \eqref{eq:laplace_coeff_err} and
  \eqref{eq:limit_gives_deriv},
  \begin{align}
    \estq(n, m) =
    \frac{r^m}{m!}
    \abs{ \frac{2n+1}{\gamma'(t_0) \sqrt{t_0^2-1}} }^m
    \frac{ |\pullb\den(t_0)| }{ |t_0 \pm \sqrt{t_0^2-1}|^{2n+1} }.
    \label{eq:am_err}
  \end{align}
  This result is essentially a generalization of
  \cite[Thm. 1]{AfKlinteberg2016quad} to curved panels.
\end{hlbox}

The estimate \eqref{eq:am_err} holds well in the limit
$n \to \infty$, but will lose accuracy with increasing $m$ for a fixed
$n$. We have in practice found that it starts losing its reliability
around $m > n/2$, so that it then is safer to trigger an upsampling in
\cref{alg:coeff}, instead of continuing to rely on the estimate.

Though the above derivation is rather straightforward, some more work
is required if we want to evaluate \eqref{eq:am_err} in practice,
since that requires finding $t_0$ and being able to evaluate
$\pullb\den(t_0)$ and $\param(t_0)$. For this, one would ideally like
to have access to analytic expressions for $\den$ and $\param$, but
what we typically have is instead the values of the functions at the
quadrature nodes. To be able to evaluate \eqref{eq:am_err} using this
pointwise data, we first need to construct continuations of $\param$
and $\pullb\den$. Then we need to solve \eqref{eq:t0_def} using the
continuation of $\param$, in order to find $t_0$ and finally evaluate
the estimate. This may sound hard, but can actually be done in an
efficient way using polynomial extrapolation.

\begin{figure}[htbp]
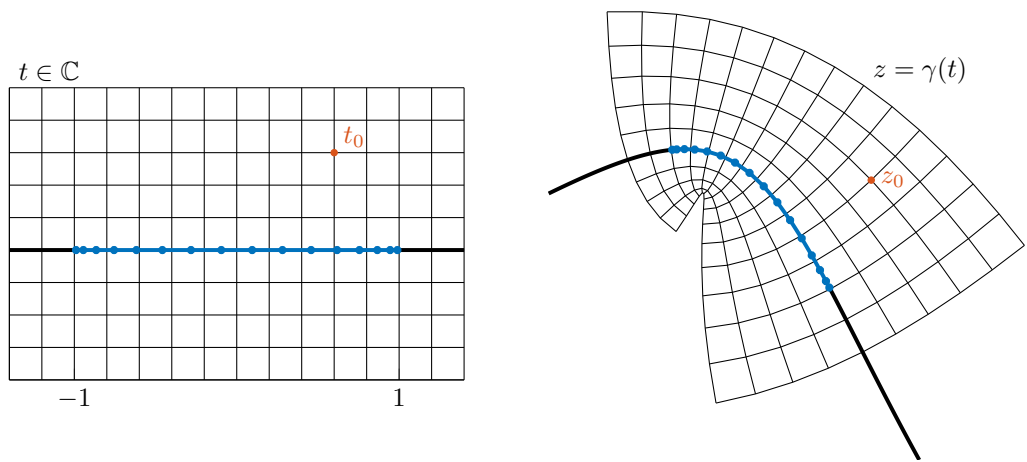

  \centering

  \setlength\fwidth{0.4\textwidth}
  \begin{tikzpicture}
    \input{single_panel_mapping_2.tikz}
    \path
    ([shift={(-5\pgflinewidth,-5\pgflinewidth)}]current bounding box.south west)
    ([shift={( 5\pgflinewidth, 1em+5\pgflinewidth)}]current bounding box.north east);
  \end{tikzpicture}
  \hspace{2em}  
  \setlength\fwidth{0.5\textwidth}
  \begin{tikzpicture}
    \input{single_panel_mapping_1.tikz}
  \end{tikzpicture}

  \caption{Illustration of how the analytic continuation of $\gamma$
    maps the vicinity of $[-1,1]$ to the space surrounding
    $\Gamma$. Here $\Gamma$ is a segment of the starfish domain seen
    in \cref{fig:dbl_lyr_err}. Polynomial interpolation of
    $\gamma$ on a 16-point panel gives a locally very accurate
    approximation of the continuation:
    $\abs{\gamma - \poly_n[\gamma]} < 10^{-7}$ on the shown grid.}
  \label{fig:grid_mapping}
\end{figure}

Beginning with the continuation of $\param$, let $t_i$ and $w_i$,
$i=1 \dots n$, be the nodes and weights of the n-point Gauss-Legendre
quadrature. We can then let $\poly_n[\gamma](t)$ be the polynomial of
degree $n-1$ that interpolates $\gamma(t)$ at the nodes $t_i$. \hl{
  For this high-order ($n=16$) interpolation to be accurate, we need
  to compute $\poly_n[\gamma]$ in a way that is both well-conditioned
  and stable \mbox{\cite[Ch. 14]{Trefethen2013}}. The distribution of
  the Legendre points $t_i$ ensures that our interpolation problem is
  well-conditioned. For stability, we use as our basis the Legendre
  polynomials $P_\ell(t)$, } %
\begin{align}
  \poly_n[\gamma](t) = \sum_{\ell=0}^{n-1} \hat\gamma_\ell P_\ell(t) .
  \label{eq:legendre_interpolant}
\end{align}
\begin{hlbox}
  These are orthogonal on $[-1,1]$, allowing us to explicitly
  compute the coefficients $\hat\gamma_\ell$ as
  \begin{align}
    \hat\gamma_\ell = \sum_{m=1}^n L_{\ell m} \gamma(t_m),
    \quad
    \ell=0, \dots, n-1,
  \end{align}
  where, for a given $n$, $L$ is a constant matrix with elements
  \begin{align}
    L_{\ell m} = \frac{2\ell+1}{2} P_\ell(t_m) w_m.
  \end{align}
\end{hlbox}
\hl{To improve convergence in the following step}, we assume that $\Gamma$
has endpoints at $-1$ and $1$; this can be achieved for any open curve
by first applying a simple scaling and rotation to both $\Gamma$ and
$z_0$. Letting the interpolant $\poly_n[\gamma](t)$ be the analytic
continuation of $\gamma(t)$, we can now find $t_0$ in
\eqref{eq:t0_def} by instead solving
\begin{align}
  \poly_n[\gamma](t_0) = z_0,
  \label{eq:t0_eq}
\end{align}
using a numerical root-finding algorithm. This does, for the purposes
of error estimation, work very well in the near neighborhood of
$\Gamma$, see example in \cref{fig:grid_mapping}. %
\hl{ The Legendre polynomials $P_\ell(t)$ are evaluated using the
  standard recurrence relation \mbox{\cite[\S14.10.3]{NIST:DLMF}},
  which holds also for complex $t$
  \mbox{\cite[\S14.21(iii)]{NIST:DLMF}}.  The polynomials grow
  exponentially in $\ell$ for $t\notin [-1,1]$, but we have not
  observed any stability issues related to this. In particular, the
  coefficients $\gamma_\ell$ of a well-resolved panel decay faster
  than the polynomials $P_\ell(t)$ grow for the range of $t$-values
  that are of interest to us. } %

\begin{figure}[htbp]
  \centering
  \includegraphics[width=.45\textwidth]{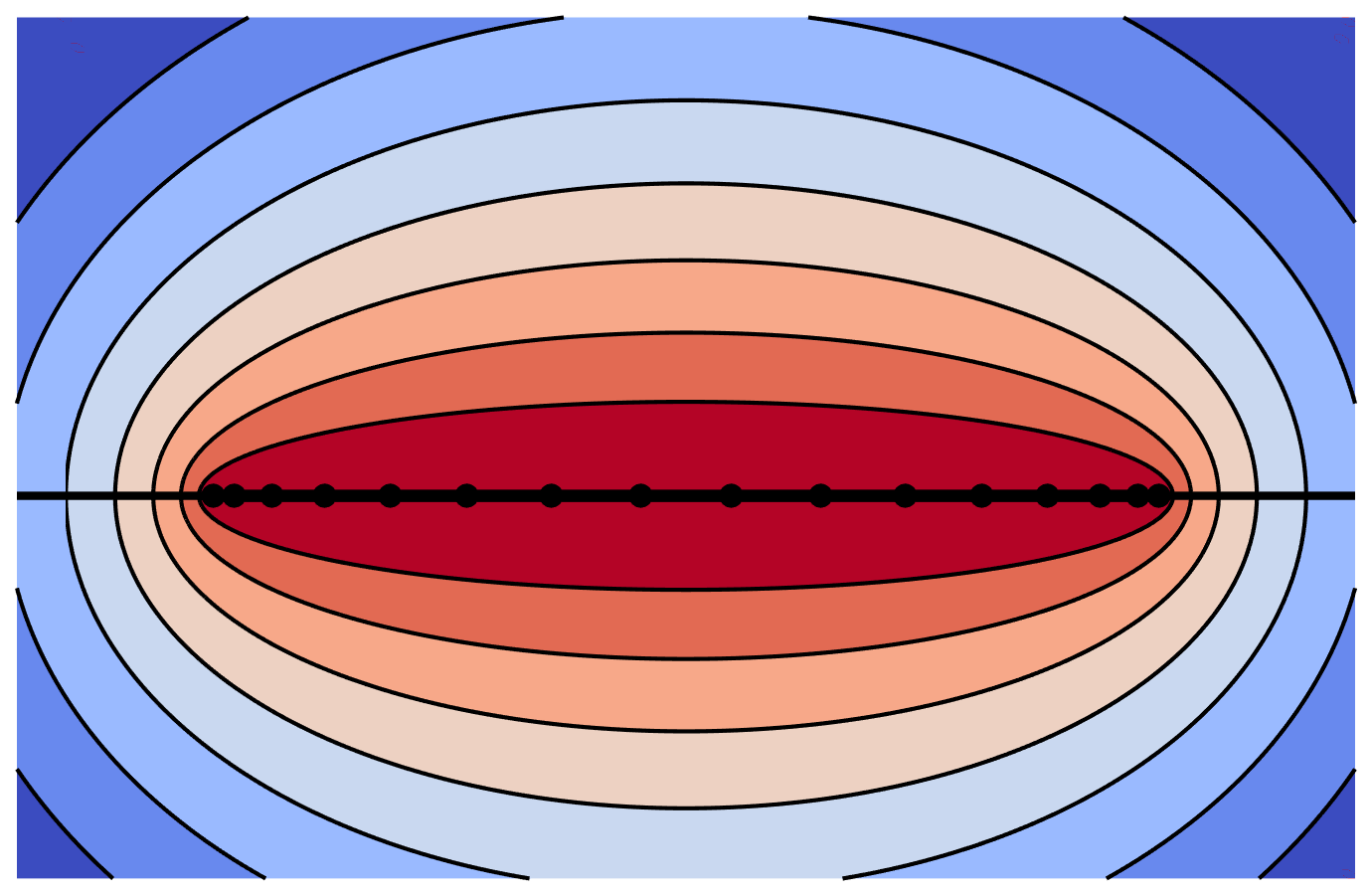}
  \hspace{2em}
  \includegraphics[width=.45\textwidth]{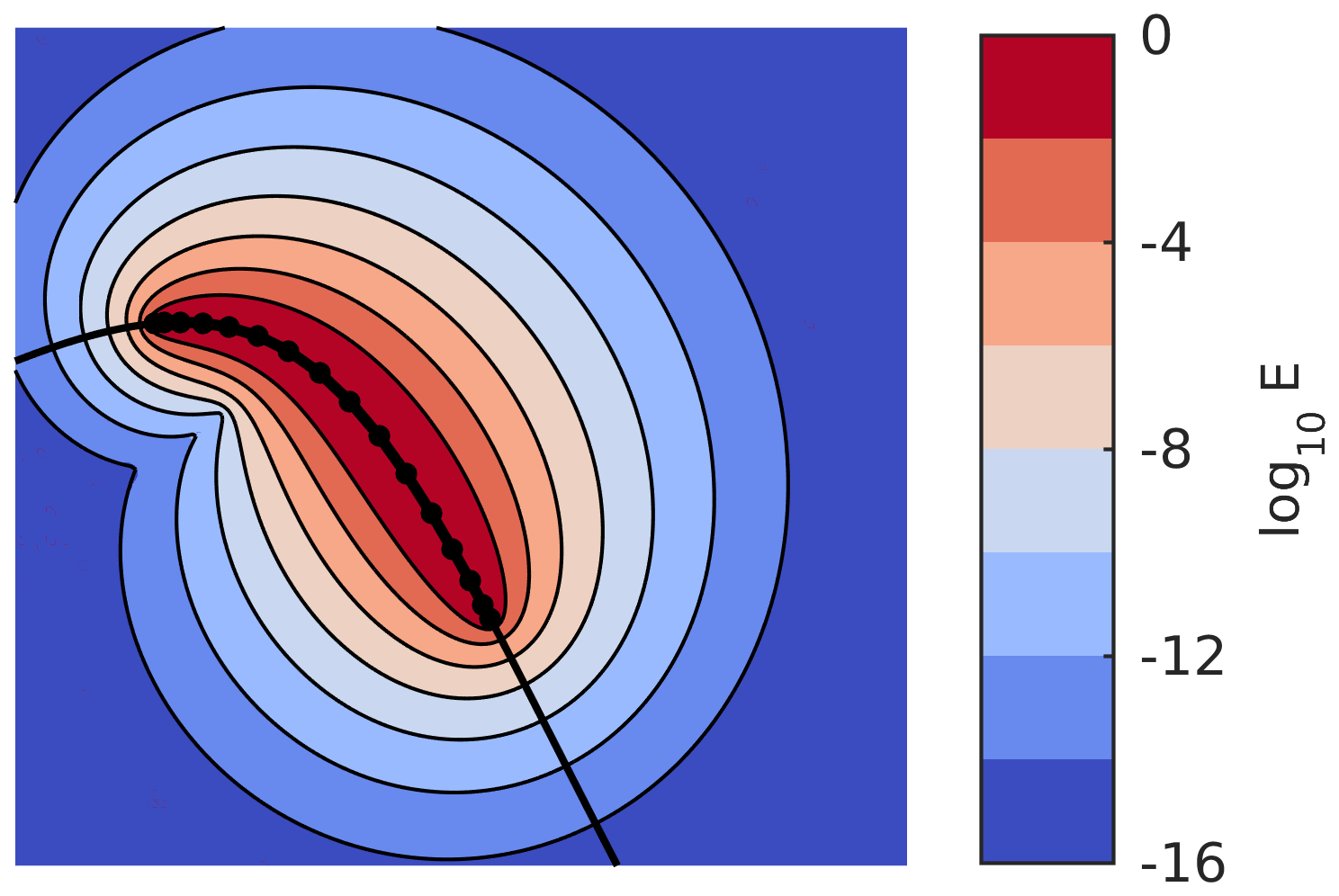}
  \caption{Quadrature errors when evaluating the integral
    \eqref{eq:dbl_lyr_coef_t} for $m=0$ and $\den=1$ on a flat and a
    curved panel. These panels correspond to those shown in
    \cref{fig:grid_mapping}.  Colored fields are actual error levels,
    contour lines are computed using the error estimate
    \eqref{eq:am_err} with $t_0$ determined by inverting
    $\poly_n[\gamma](t)$. The ellipses to the left are the Bernstein
    ellipses.}
  \label{fig:panel_err}
\end{figure}

Equation \eqref{eq:t0_eq} can be solved efficiently using Newton's
method, but the choice of starting guess can be important, especially
near concave regions of the curve (e.g. below the curve of
\cref{fig:grid_mapping}). In such regions the inverse mapping
$t_0 = \gamma^{-1}(z_0)$ is no longer single-valued, and for our
estimate to be accurate we want the solution $t_0$ that predicts the
largest error, when the estimate is evaluated at that point. This
``worst solution'' can often be found by solving with $t_0 = \pm 1$ as
starting guesses and comparing the results; this strategy was used for
the results shown to the right in \cref{fig:panel_err}. As an
alternative, one can use $t_0=z_0$ as as starting guess (given the
above assumptions on the endpoints of $\Gamma$). This simpler strategy
works well for practical purposes, and is what we use throughout the
remainder of this paper.


\begin{figure}[htbp]
  \begin{hlbox}
  \centering
  \begin{minipage}{0.45\textwidth}  
    \includegraphics[width=\textwidth]{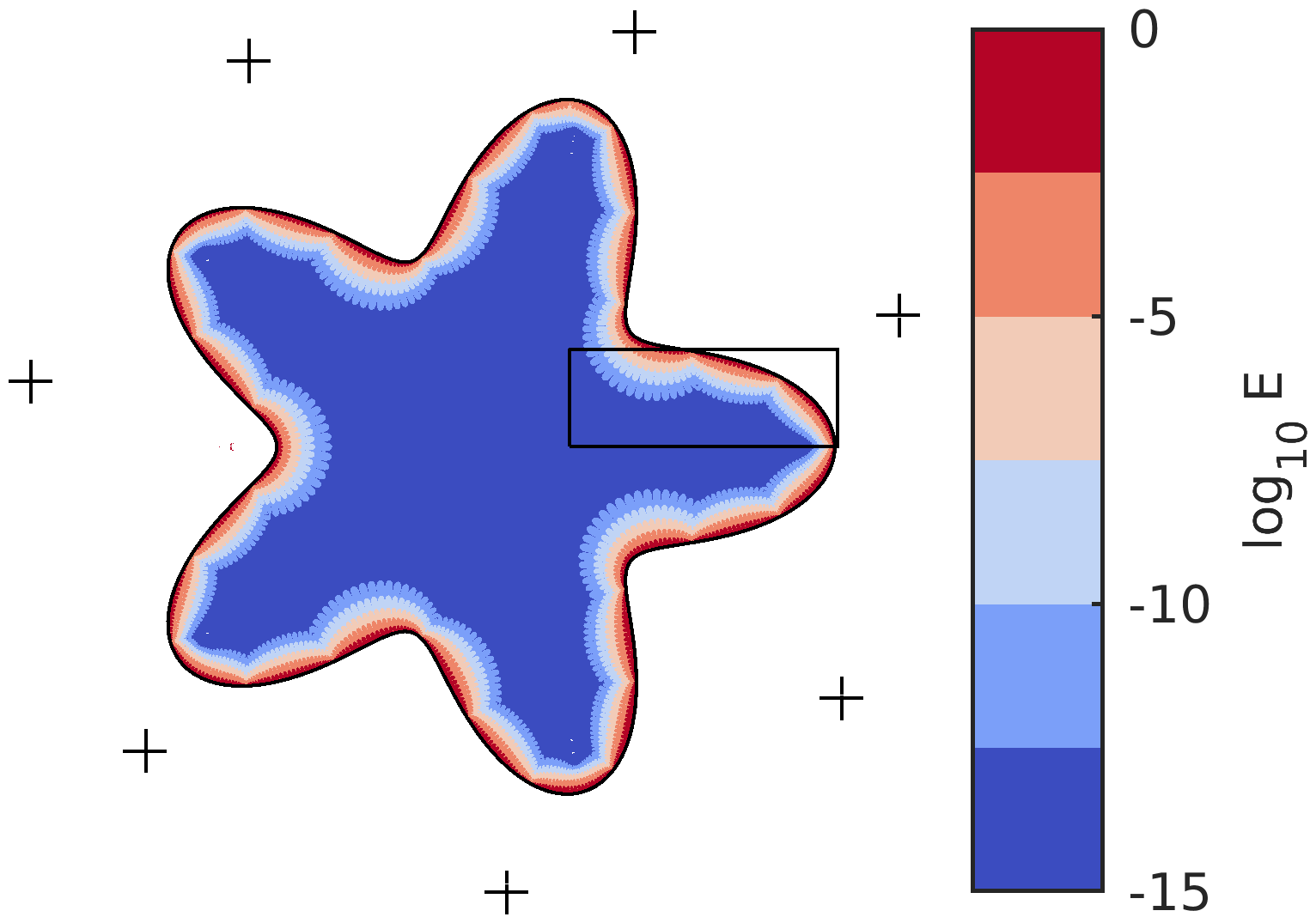}
  \end{minipage}
  \begin{minipage}{0.54\textwidth}
    \includegraphics[width=\textwidth]{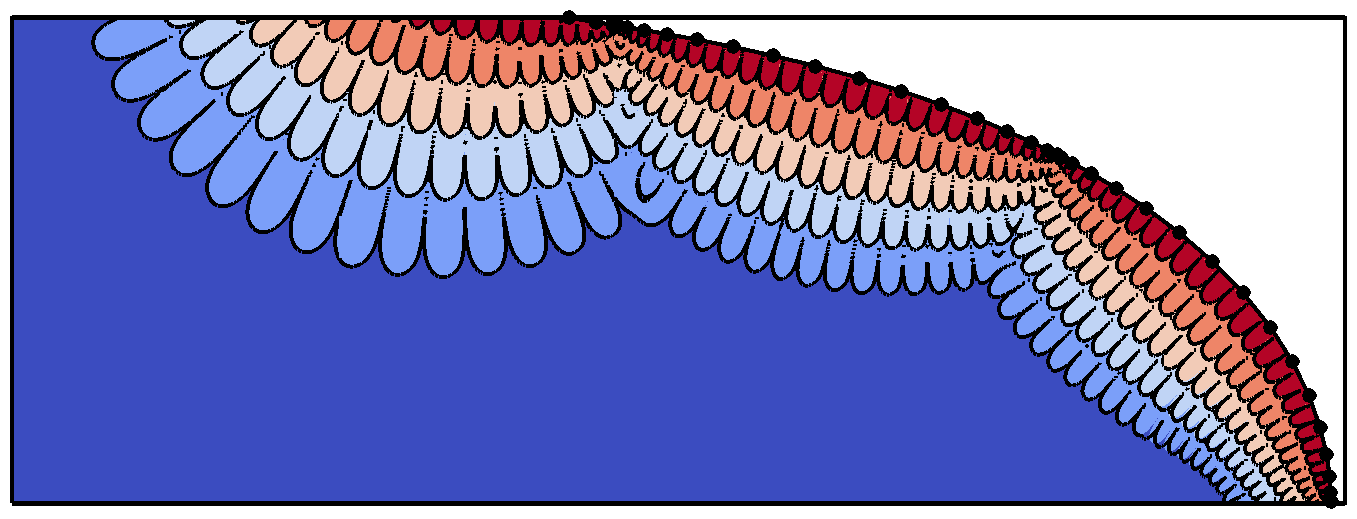} \\
    \includegraphics[width=\textwidth]{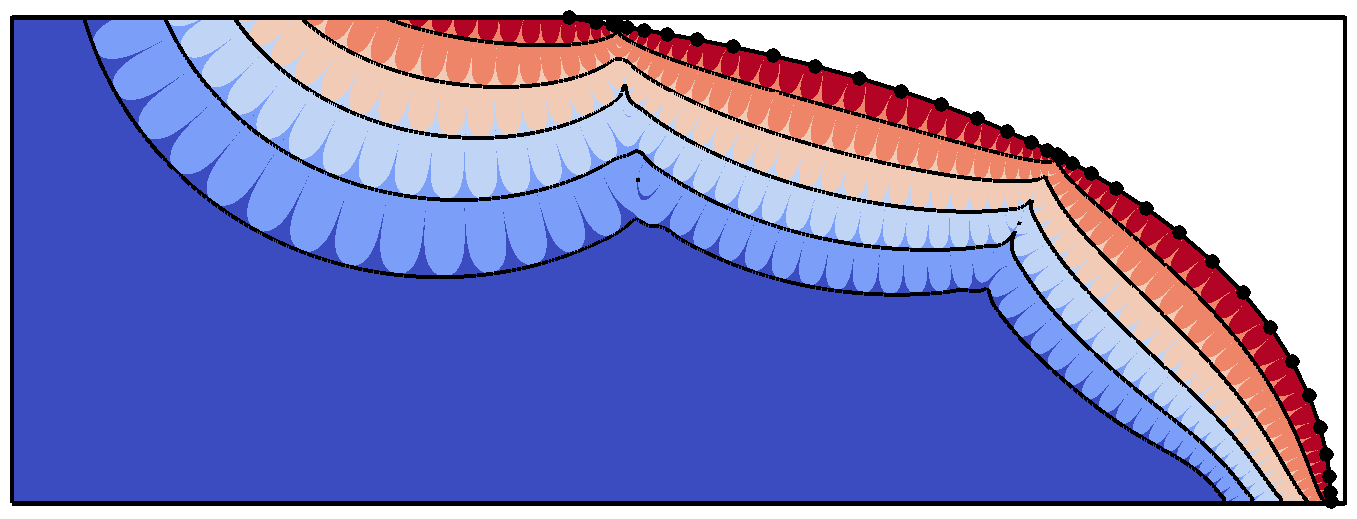}    
  \end{minipage}
  \caption{Error curves when evaluating the Laplace double layer
    potential directly using 27 panels of equal arc length, with 16
    points on each panel. Colored fields represent the error compared
    to the exact solution. Black contour lines in the top right are computed
    using the estimate \eqref{eq:dbl_lyr_err} with
    \eqref{eq:den_t0_extrap}, and coincide almost perfectly with the
    actual error. Black contour lines in the bottom right are computed
    without the $\Im[\cdot]$ in the estimate \eqref{eq:dbl_lyr_err},
    which produces smooth contours that bound the error curves rather
    than lie on top of them. This would be enough in most practical
    applications.}
  \label{fig:dbl_lyr_err}
\end{hlbox}

\end{figure}

Once $t_0$ is found, we are able to evaluate the coefficient error
estimate as given in \eqref{eq:am_err}. The value of $\param'(t_0)$
can be found by differentiation of the interpolant,
\begin{align}
  \param'(t_0) = P_n'[\gamma](t_0) .
\end{align}
We can also evaluate $\pullb\den(t_0)$ through an interpolating polynomial,
\begin{align}
  \pullb\den(t_0) = P_n[\pullb\den](t_0).
  \label{eq:den_t0_extrap}
\end{align}
This works very well if $\den$ is well-resolved on $\Gamma$, and can
be used to obtain the fine-scale correspondence between error and
estimate seen in \cref{fig:dbl_lyr_err}. A less expensive
alternative is to use the max norm of $\pullb\den$ on the
interval/panel,
\begin{align}
  \pullb\den(t_0) \approx \| \pullb\den(t) \|_{L^\infty(-1, 1)} 
  = \|   \den(z) \|_{L^\infty(\Gamma)}.
  \label{eq:den_t0_max}
\end{align}
\hl{This works well in practice, as we mainly need to get the order of
  magnitude right, and also appears to be slightly more robust
  whenever $\den$ is not fully resolved on $\Gamma$, especially for
  Helmholtz (see
  \mbox{\cref{fig:helmholtz_err,fig:helmholtz_err_close}}). We use
  this approach in the results of section
  \mbox{\ref{sec:numer-exper}}}.

The above results can also be used for estimating the nearly singular
quadrature error when evaluating the double layer potential
\eqref{eq:laplace_dbl_lyr} near $\Gamma$, and hence to determine when
AQBX needs to be used. To that end, we simply use the observation that
$u(z_0) = \Re \coeff_0$. Denoting by $\tilde u(z_0)$ the potential
evaluated using direct quadrature, we thus have from
\eqref{eq:laplace_coeff_err} that
\begin{align}
  \abs{u(z_0) - \tilde u(z_0)} \approx
  \frac{1}{2\pi} \abs{\Im \left[\pullb\den(t_0) k_n(t_0)\right]} .
  \label{eq:dbl_lyr_err}
\end{align}
This estimate, when evaluated using the above procedure, is very
accurate. As an example, in  \cref{fig:dbl_lyr_err}, we consider
the Laplace Dirichlet problem in a starfish domain. The density is
computed by solving \eqref{eq:laplace_inteq} using the Nystr\"om
method and a right-hand side given by a collection of point sources
(marked by + in the figure). The solution is then evaluated inside the
domain directly using the composite Gauss-Legendre quadrature, and
compared to the reference solution given by the boundary
condition. Comparing these results to those in
\cite[fig. 7]{AfKlinteberg2016quad}, which used a flat panel
approximation, shows the importance of taking into the account the
inverse mapping of the target point.

\subsection{Helmholtz coefficient error}

We now turn to AQBX for the Helmholtz combined field potential
(sec. \ref{sec:helmholtz-equation}). When evaluating the expansion
coefficients, the Hankel functions in the integrands have a
singularity as $r_w \to 0$, which to leading order behaves like
\cite[\S10.4,\S10.8]{NIST:DLMF}
\begin{align}
  \hankfk_{\pm l}(kr_w) = -(\pm 1)^l \frac{2^l(l-1)!}{\pi}(kr_w)^{-l} 
  + \ordo\left((kr_w)^{-(l-2)}\right).
  \label{eq:hankfk_exp}
\end{align}
The quadrature error due to near singularity is dominated by that from
the highest-order pole, so for the purposes of error estimation it is
suitable to only keep the highest-order Hankel function in the
expression for \eqref{eq:addA_helm}, and to approximate that Hankel
function using only the first term of \eqref{eq:hankfk_exp}. We can
then approximate \eqref{eq:addAS} as
\begin{align}
  \addAS_m(w, z_0) &\approx 0,
\end{align}
and \eqref{eq:addAD} as
\begin{align}
  \begin{split}    
    \addAD_0(w, z_0) &= -\frac{ik}{4} \hankfk_1(kr_w) r_w \Re \left[ \frac{n_w}{w-z_0} \right]
    \\
    &\approx 
    \frac{1}{2\pi} \Re \left[ \frac{n_w}{w-z_0}\right]
  \end{split}
      \label{eq:d0_approx}
\end{align}
and, for $|m|>0$ and using \eqref{eq:exppot},
\begin{align}
  \addAD_m(w, z_0) &\approx \frac{k}{8}\hankfk_{|m|+1}(kr_w) e^{-i(|m|+1)\theta_w}n_w \\
                   &\approx \frac{k}{8\pi} 2^{|m|+1}|m|! (kr_w)^{-(|m|+1)} e^{-i(|m|+1)\theta_w}n_w \\
                   &= \frac{2^{|m|}|m|!}{4\pi k^{|m|}} \frac{n_w}{(w-z_0)^{|m|+1}}.
\end{align}
\begin{hlbox}
  Inserted into \eqref{eq:addA0_helm} and \eqref{eq:addA_helm}, this gives
  \begin{align}
    \addA_0^r(w, z_0) &\approx \frac{1}{2\pi} 
                        \Re \left[ \frac{n_w}{w-z_0}\right],\\
  \addA_m^r(w, z_0) &\approx 
                      \sqrt{2} \left(1, 1 \right) \frac{r^{m}}{4\pi}
                      \frac{n_w}{(w-z_0)^{m+1}}, \quad m>0 .
                      \label{eq:addAm_approx}
\end{align}
\end{hlbox}
Note that we in the above series of simplifications
(\ref{eq:d0_approx}--\ref{eq:addAm_approx}) have removed constant
factors which are irrelevant to the magnitude of the error (i.e. $-1$
and $i$). Estimating the error of these simplified integrands on a
panel with parametrization $\param$, it follows that
\begin{align}
  \label{eq:helm_coeff_err_zero}
  |\coeff_0-\tilde\coeff_0| &\approx \abs{\Im \oprem\left[ 
                              \frac{\pullb\den \param'}{ 2\pi (\param-z_0)}
                              \right]}, \\
  \label{eq:helm_coeff_err}
  |\coeff_m-\tilde\coeff_m| &\approx \abs{\oprem\left[ \frac{r^m \pullb\den \param'}{2\pi (\param-z_0)^{m+1}} \right]}, \quad m>0.
\end{align}
Note that \eqref{eq:helm_coeff_err_zero} is bounded by
\eqref{eq:helm_coeff_err} with $m=0$, so it is sufficient to use only
\eqref{eq:helm_coeff_err}, if we can accept being on the conservative
side for $m=0$. Also note that the integrand in
\eqref{eq:helm_coeff_err} is exactly the same as in the error estimate
for the Laplace double layer potential
\eqref{eq:laplace_coeff_err}. The conclusion is that, remarkably, the
coefficient error estimate $\estq$ for Helmholtz is identical to that
previously derived for Laplace \eqref{eq:am_err},
\begin{align}
  \estq(n, m) =
  \frac{r^m \abs{\pullb\den(t_0)}}
  {2\pi m! \abs{\gamma'(t_0)}^{m}}
  \abs{k_n^{(m)}(t_0)}.
  \label{eq:am_err_helmholtz}
\end{align}
This correspondence between Laplace and Helmholtz also holds for the
quadrature error when evaluating the underlying potential itself. This
can be seen from \eqref{eq:helm_coeff_err_zero} and the observation
from \eqref{eq:addADzero} that \hl{$a_0=u(z_0)$},
\begin{align}
  \abs{u(z_0) - \tilde u(z_0)} &\approx
  \frac{1}{2\pi} \abs{\Im \left[\pullb\den(t_0) k_n(t_0)\right]} 
  \label{eq:helmholtz_err_im} \\
  &\le \frac{1}{2\pi} \abs{\pullb\den(t_0) k_n(t_0)}.
  \label{eq:helmholtz_err}
\end{align}
This estimate has been derived by simplifying the Helmholtz double
layer kernel \eqref{eq:d0_approx}, and noting that the singularity to
leading order is identical to that of the Laplace double
layer. However, in experiments we can observe that the small-scale
oscillations predicted by \eqref{eq:helmholtz_err_im} (and which are
clearly noticeable in the Laplace case) only appear for small
wavenumbers, and it is therefore generally better to use
\eqref{eq:helmholtz_err}. \hl{ The reason for the disappearance of the
  oscillations is unknown to us, though an interesting feature is that
  they only seem to disappear when $\den$ is the solution to
  \mbox{\eqref{eq:helmholtz_inteq}}. Setting $\den=1$ and computing
  the error by comparing to a finer grid produces the oscillations
  also for large wavenumbers.}

\begin{figure}[htbp]
  \centering
  \begin{hlbox}
    \includegraphics[width=.36\textwidth]{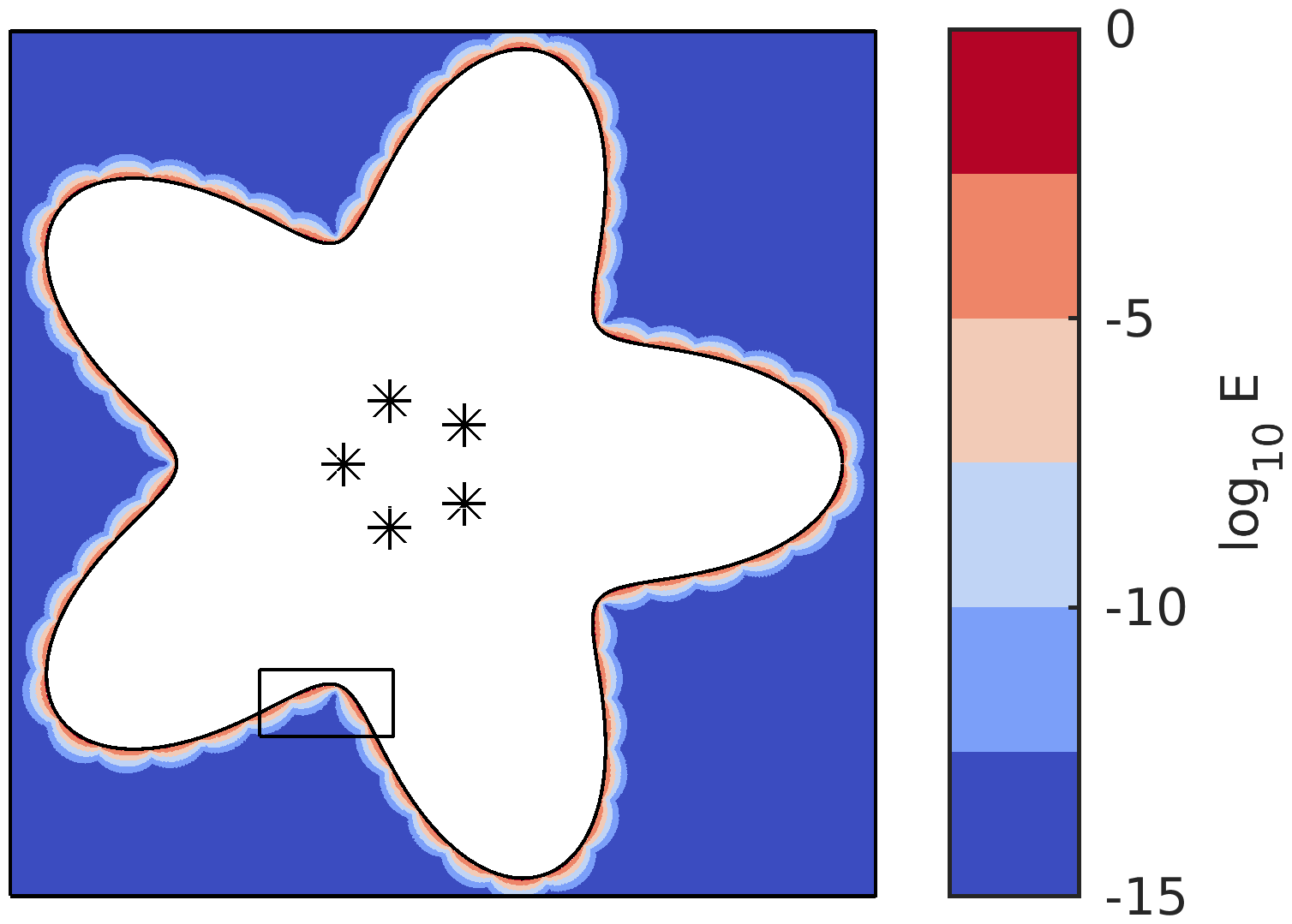}
    \hspace{0.5cm}
    \includegraphics[width=.5\textwidth]{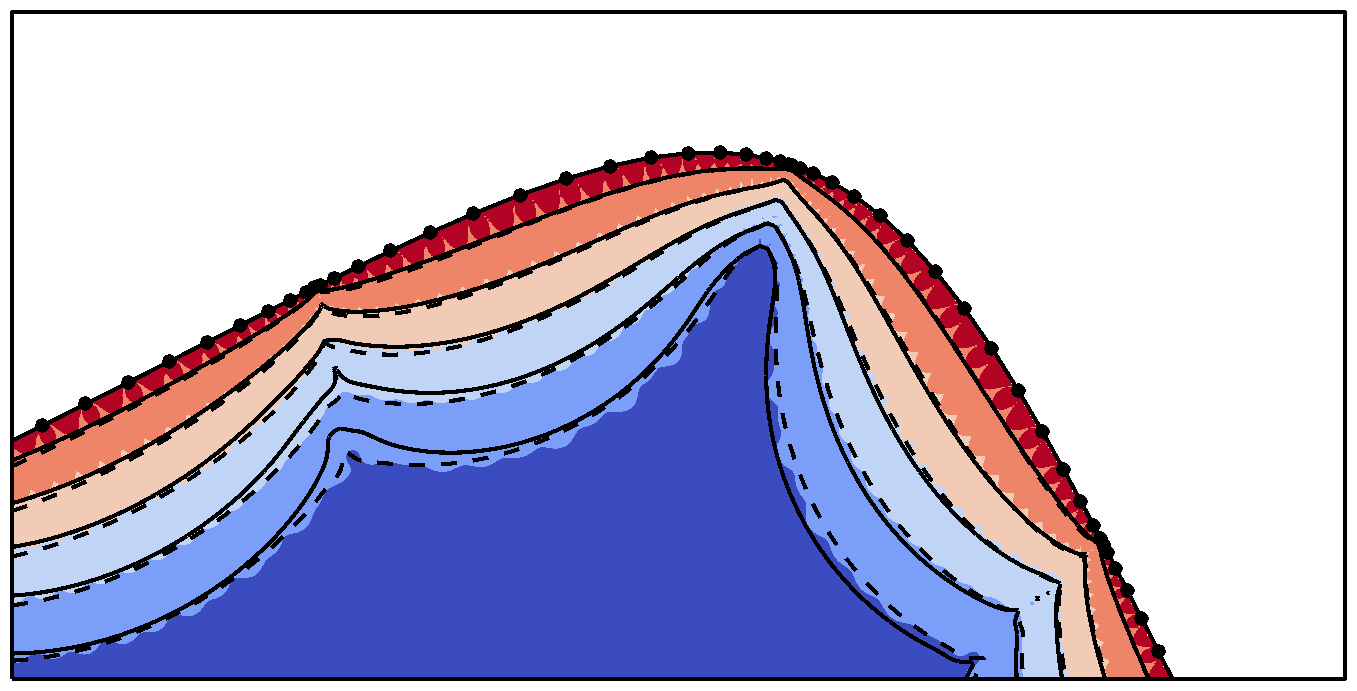}
    \caption{Error curves when evaluating the Helmholtz combined field
      potential directly using 60 panels of equal arc length, with 16
      points on each panel and the point sources located on a circle of
      radius 0.2. The wavenumber is set to $k=4/h=26.6$, $h$ being the length
      of the panels. Colored fields represent the error compared to the
      exact solution, black contours are computed using the estimate
      \eqref{eq:helmholtz_err}. The solid lines correspond to
      \eqref{eq:helmholtz_err} with $\den(t_0)$ evaluated using
      \eqref{eq:den_t0_extrap}, while the dashed lines use the
      approximation \eqref{eq:den_t0_max}.}
    \label{fig:helmholtz_err}
  \end{hlbox}
\end{figure}

\begin{figure}[htbp]
  \centering
  \begin{hlbox}
    \includegraphics[width=.36\textwidth]{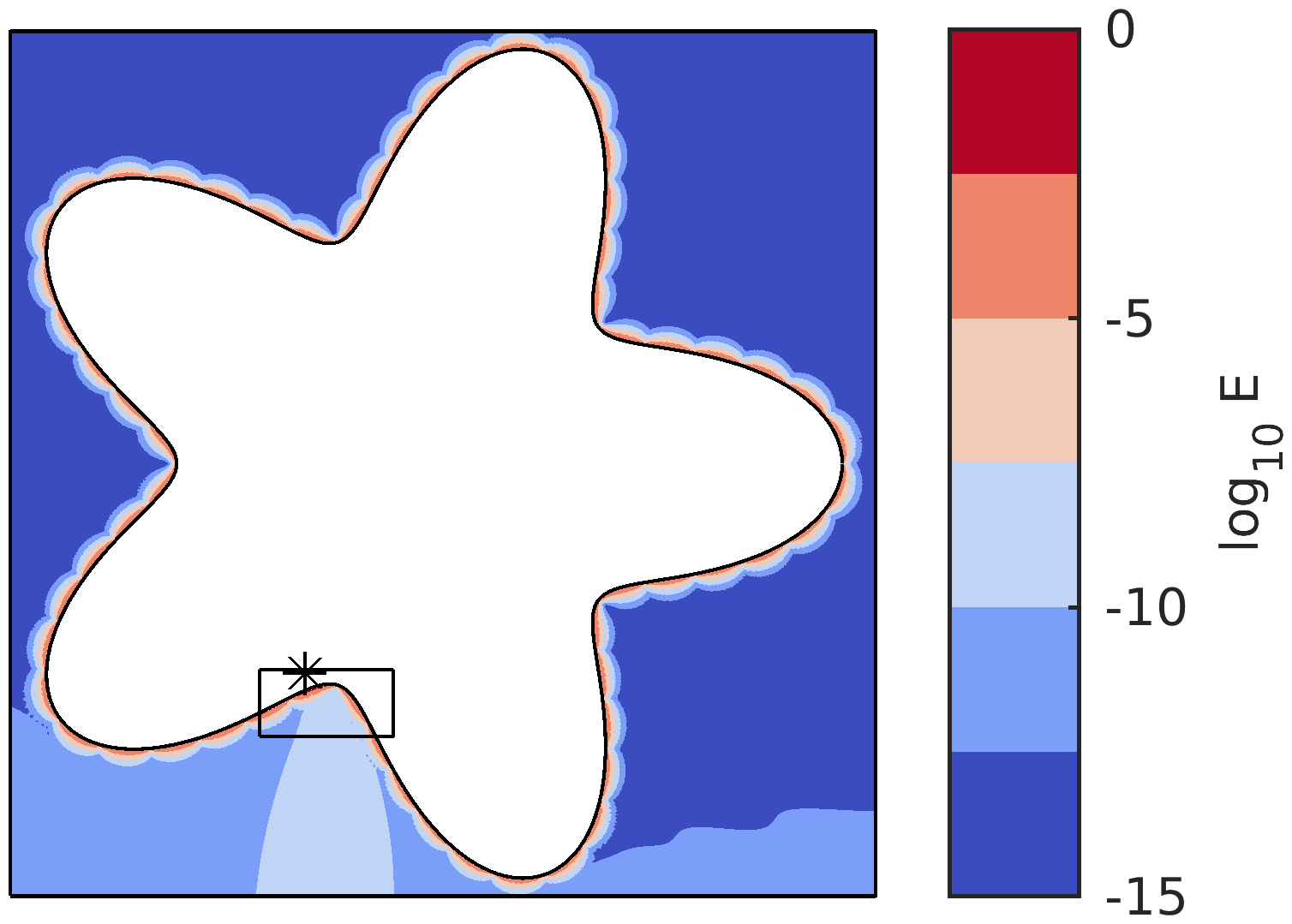}
    \hspace{0.5cm}    
    \includegraphics[width=.5\textwidth]{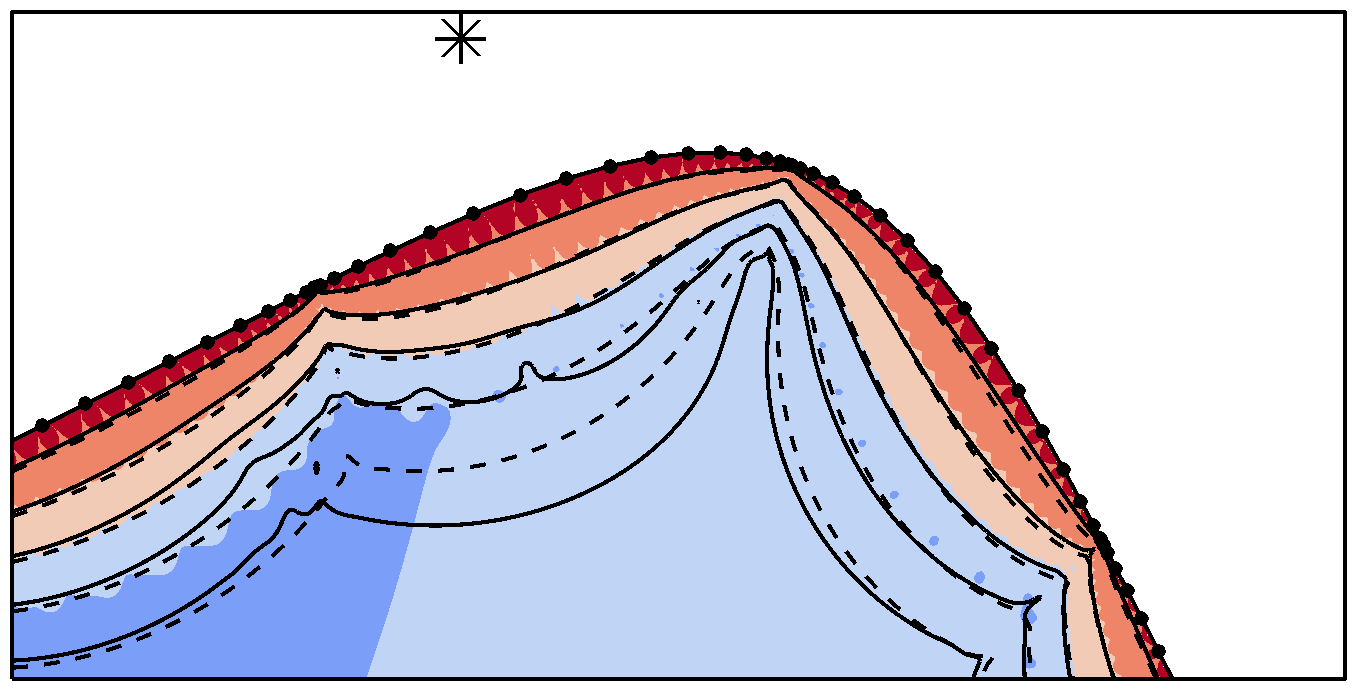}
    \caption{Here the discretization and visualization are the same as
      in \cref{fig:helmholtz_err}, but the source is now located at a
      distance $h/3$ away from the panel. The error is still well
      estimated close to the boundary, where it is dominated by the
      nearly singular quadrature error, though evaluating $\den(t_0)$
      using the polynomial extrapolation \eqref{eq:den_t0_extrap}
      appears unstable. Further away from the boundary the error is
      dominated by a lack of resolution, which our estimate will not
      capture.}
    \label{fig:helmholtz_err_close}
  \end{hlbox}
\end{figure}

As a demonstration, in \cref{fig:helmholtz_err} we repeat the
experiment of \cref{fig:dbl_lyr_err}, but for the Helmholtz exterior
Dirichlet problem. We set a number of point sources (marked $*$)
inside a starfish domain and solve the integral equation
\eqref{eq:helmholtz_inteq} using the discretization scheme of
\cite{Helsing2015}. The correspondence between error and estimate is
still very good, though not as excellent as in the Laplace case.
\hl{\mbox{\cref{fig:helmholtz_err_close}} shows that the estimate
  appears to work also for a source point very close to the boundary,
  when the accuracy of the solution $\den$ has started to
  deteriorate. Here, when $\den$ is no longer well-resolved, one can
  clearly see that evaluating $\den(t_0)$ using the panel max norm
  \mbox{\eqref{eq:den_t0_max}} is more stable than using polynomial
  extrapolation \mbox{\eqref{eq:den_t0_extrap}}. Using
  \mbox{\eqref{eq:den_t0_max}} results in a useful error estimate
  close the boundary, where the nearly singular quadrature error
  dominates. Further away the error is due to $\den$ not being
  accurate, something which our estimate does not take into account.}

It is worth noting that the error estimate \eqref{eq:helmholtz_err} is
independent of the wavenumber $k$. This might come as a surprise, as
one usually needs to increase the grid resolution with increasing
wavenumber. However, here we only take into account the nearly
singular quadrature error, under the assumption that far-field
interactions are well-resolved. The result simply reflects the fact
that the singularity in the kernel is independent of $k$.

\section{Local AQBX}
\label{sec:local-aqbx}
\newcommand{\near}{{\mathcal N}}
\newcommand{\far}{{\mathcal F}}

Since QBX is a special quadrature scheme for target points that are
close to or on the boundary $\bdry$, it makes sense to only use QBX
for those parts of the boundary that are close to a given target
point. This is known as ``local QBX'' \cite{Rachh2015a} (as opposed to
``global QBX''), and can be particularly straightforward to combine
with a fast method. For panel-based quadrature on a simple curve this
is easy to implement; only the panels that are near a given expansion
center are used in the local expansion. Selecting panels to include
can be done using an error estimate of the layer potential, such as
\eqref{eq:dbl_lyr_err} in combination with a tolerance, or by simply
including a fixed number of neighboring panels (this works well if all
panels are of equal length). When evaluating the potential, the
contribution from the near panels is computed through the local
expansion, while the contribution from the remaining panels is
computed directly using the underlying Gauss-Legendre quadrature.

Let the boundary be composed of a set of panels $\Gamma_i$,
\begin{align}
  \bdry = \bigcup_i \Gamma_i .
\end{align}
We can then denote by $\near$ the near panels that are included in the
local expansion, and by $\far$ the far panels that are evaluated
directly. \hl{ Note that this division must be made such that the
  panels in $\far$ are well-separated from all target points at which
  the expansion will be evaluated. }  If we write the layer potential
as
\begin{align}
  u^\bdry(z) = \int_\bdry G(z, w) \den(w) \dif s_w,
\end{align}
then the numerical approximation of $u$ using local QBX can be written
as
\begin{align}
  \tilde u^\bdry(z) = 
  u_{\text{QBX}}^\near(z) + u_{\text{direct}}^\far(z) .
\end{align}
To combine this with a fast method that directly evaluates the
interactions between all source points (such as the FMM), one can
simply subtract the direct contribution from the near panels,
\begin{align}
  \tilde u^\bdry(z) = 
  u_{\text{fast}}^\bdry(z) - u_{\text{direct}}^\near(z) + u_{\text{QBX}}^\near(z) .  
\end{align}
The last two terms in this expression can together be viewed as a
correction term to the direct quadrature. Since this correction only
has local support, the local QBX scheme is FMM compatible, in the
sense introduced in \cite{Hao2014}. This ``black box'' way of using
the fast method could potentially introduce cancellation errors, when
the direct contribution from the near panels is added and subtracted,
though we have not observed any such problems in practice. The
alternative is to modify the FMM to ignore those local contributions
in the first place, which is non-trivial for complex geometries.

A subtle feature of local QBX is that the width of the segment $\near$
affects the convergence rate of the local expansion. This is due to
artificially induced endpoint singularities at the interfaces between
$\near$ and $\far$, and has been discussed to some extent in both
\cite{Barnett2014} and \cite{Rachh2015a}. The choice of width of
$\near$ is therefore a balance; widening $\near$ means more points
contributing to each expansion, while narrowing $\near$ gives a slower
convergence of the expansion. Actually, not only is the convergence
slower, the exponential convergence of the truncation error
\eqref{eq:eT_bound} also tends to be less regular. This in turn makes it
harder for AQBX to correctly determine when to terminate. We have
found that a good balance is struck by using the five panels that are
closest to the expansion center.

\section{Numerical experiments}
\label{sec:numer-exper}

We have implemented the above algorithms in Matlab, and used the FMM
as implemented in FMMLIB2D \cite{fmmlib2d} for fast far-field
evaluations. Timings will not be reported here, as our code is a proof
of concept rather than a production implementation. We will in the
following numerical experiments only report on the Helmholtz problem,
as that is the more challenging one. Carrying out the same experiments
for the Laplace problem just results in similar, though slightly
better, results.

\begin{figure}[htbp]
  \begin{hlbox}
    \centering
    \hspace{1em}
    \includegraphics[height=.3\textwidth]{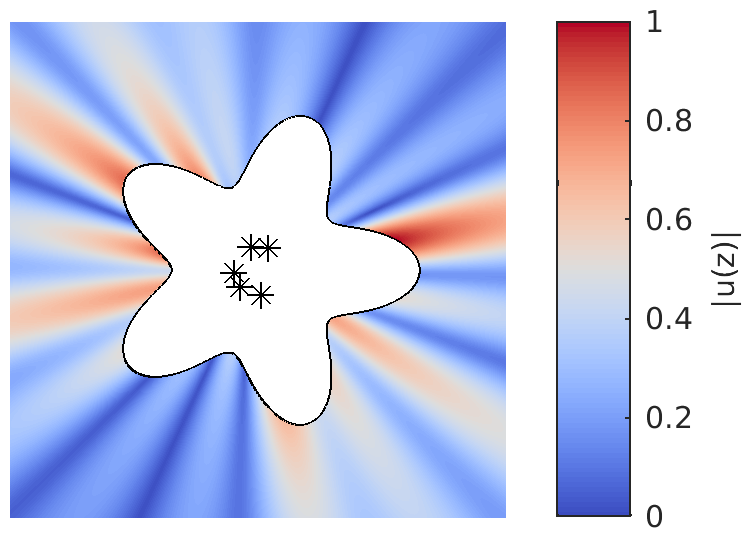}
    \hfill
    \includegraphics[height=.3\textwidth]{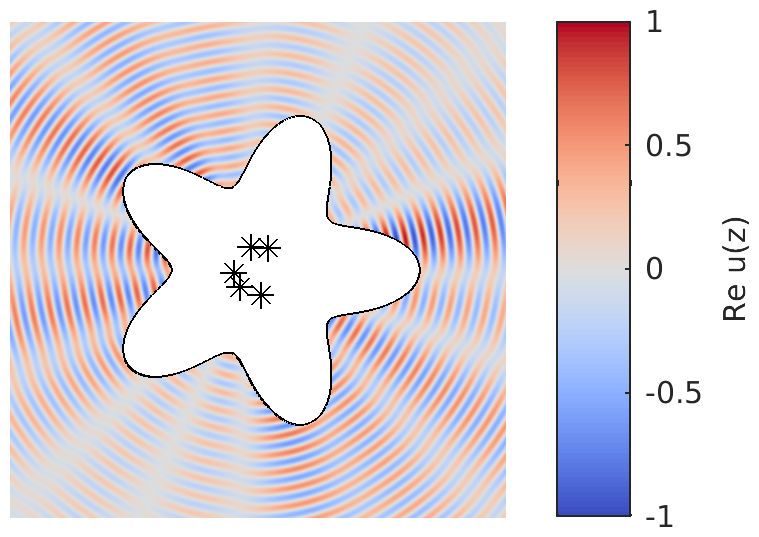}
    \hspace{0.5em}  
    \caption{Solution to Helmholtz equation ($\wnum=44.36$) given by five point sources
      located inside a starfish domain.}
    \label{fig:helmholtz_sol}
  \end{hlbox}
\end{figure}

\hl{%
  In our numerical experiments we will mainly use AQBX as a
  post-processing tool, meaning that we use it to evaluate the layer
  potential given by a known density. This allows us to study the
  behavior of AQBX in isolation, without taking into account the
  method used to compute the density. However, in section
  \mbox{\ref{sec:solv-integr-equat}} we show that AQBX can be used also for
  solving the integral equation.
  }

For our experiments we set up the reference problem shown in
\cref{fig:helmholtz_sol}: The Helmholtz problem in the domain exterior
to the starfish curve $\gamma(t) = (1+0.3\cos(10 \pi t))e^{-2\pi i
t}$, $t\in[0,1]$, with Dirichlet boundary condition given by the
potential from five point sources in the interior domain. \hl{ The
point sources are randomly positioned on a circle with radius 0.2
centered at the origin, with strengths that are randomly drawn and
then normalized such that $\norm{u}_\infty=1$ on $\bdry$. This way all
errors reported below are both relative and absolute.}

\hl{We discretize the boundary using 200 Gauss-Legendre panels of
order 16 and equal arc length $h$, and position one expansion center
at a distance $r=h/4$ in the normal direction from each point on the
boundary. The wavenumber is set to $\wnum=2/h=44.36$. The density
$\den$ is computed by solving the integral equation
\mbox{\eqref{eq:helmholtz_inteq}} using the Nystr\"om method of
\mbox{\cite{Helsing2015}}, which gives us a solution that has a
relative error of approximately $10^{-14}$ when using direct
quadrature away from the boundary (measured on a circle of radius
2). Increasing the number of panels or decreasing $\wnum$ does not
significantly improve the accuracy of the solution, nor the accuracy
of the QBX evaluation.}

\hl{Once we have computed $\den$, we can evaluate the layer potential
using AQBX, and compare the result to the exact solution given by the
potential from the five point sources. The error in AQBX has in our
tests always been largest when evaluating the layer potential on the
boundary (where the integral is singular), so we mainly report the
errors as measured there.}

\subsection{Performance of the algorithm}
\label{sec:perf-algor}

\begin{figure}[htbp]
  \centering
  \tikzset{font=\footnotesize} 
  \tikzset{mark size=1/2}
  \includegraphics[width=.38\textwidth, height=.41\textwidth]{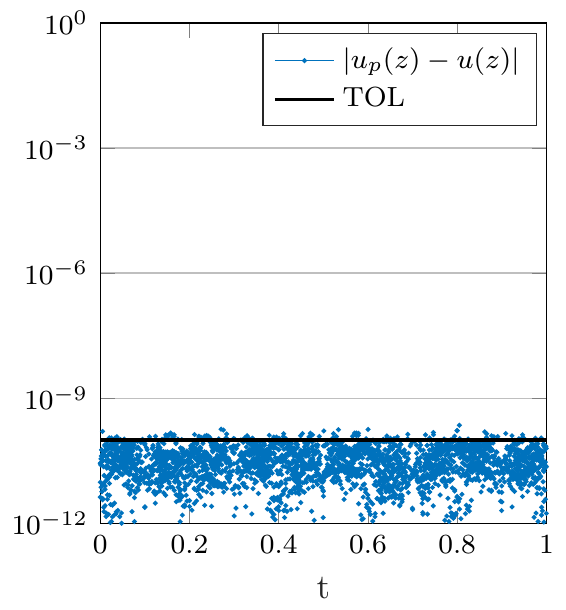}
  \tikzset{mark size=2}
  \includegraphics[width=.6\textwidth, height=.4\textwidth]{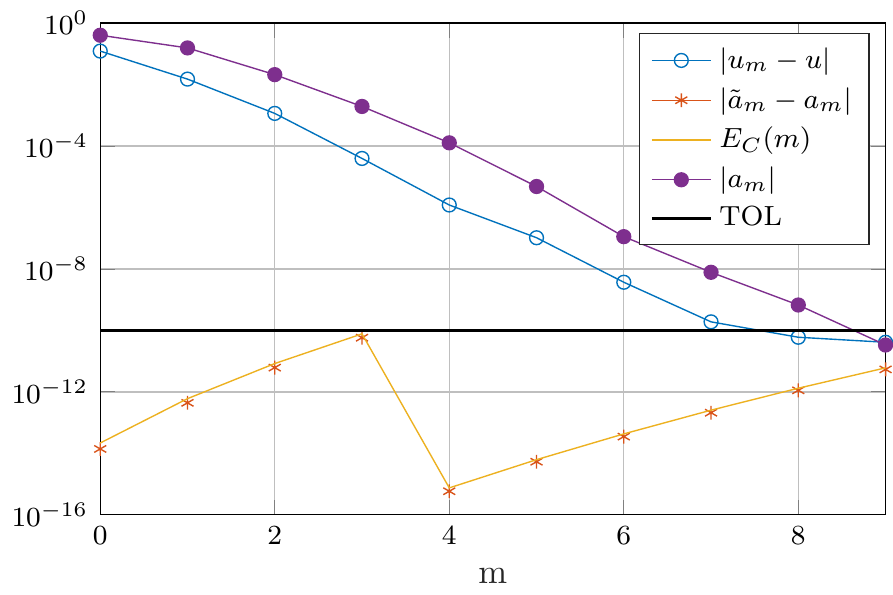}
  \caption{Results when evaluating the solution to Helmholtz equation
    using AQBX with tolerance set to $10^{-10}$, marked as thick black
    line. \emph{Left:} Distribution of error along $\bdry$, which
    shows that the error stays close to the tolerance, though it is
    not strictly met. \emph{Right:} Example showing the behavior of
    AQBX at a single expansion center, when evaluating at the closest
    boundary point. The error decays exponentially with expansion
    order, at the same rate as the coefficients $a_m$.  At the same
    time the coefficient error $|a_m - \tilde a_m|$ is growing, but is
    well estimated by the estimate $\estq(m)$ \eqref{eq:helm_coeff_err}.
    Note the jump in coefficient error between $m=3$ and $m=4$, where
    the grid is upsampled to maintain the error below tolerance.}
  \label{fig:error}
\end{figure}

\begin{figure}[htbp]
  \begin{hlbox}
    \centering
    \begin{minipage}{0.28\textwidth}
      \includegraphics[width=\textwidth]{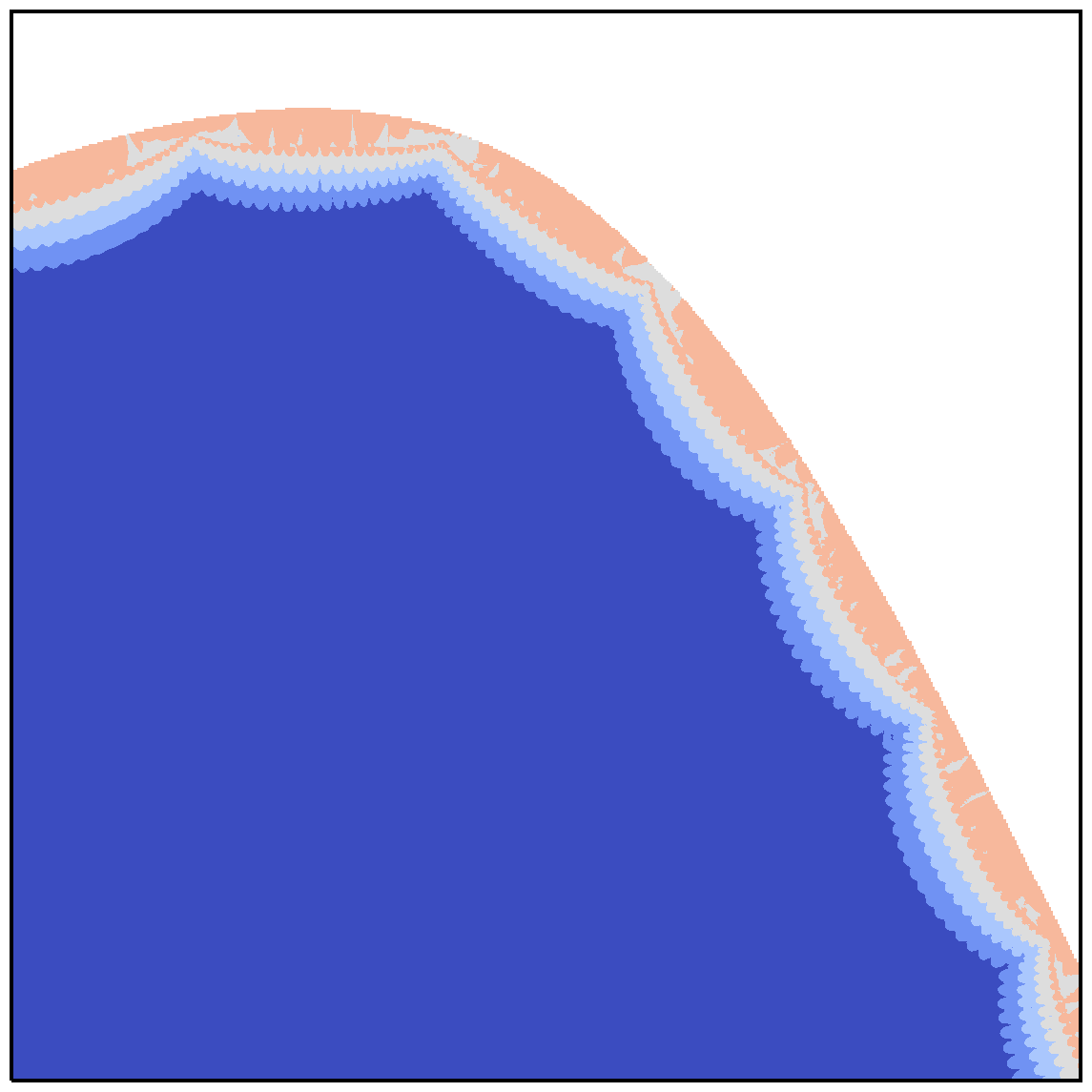}
    \end{minipage}
    \begin{minipage}{0.28\linewidth}
      \includegraphics[width=\textwidth]{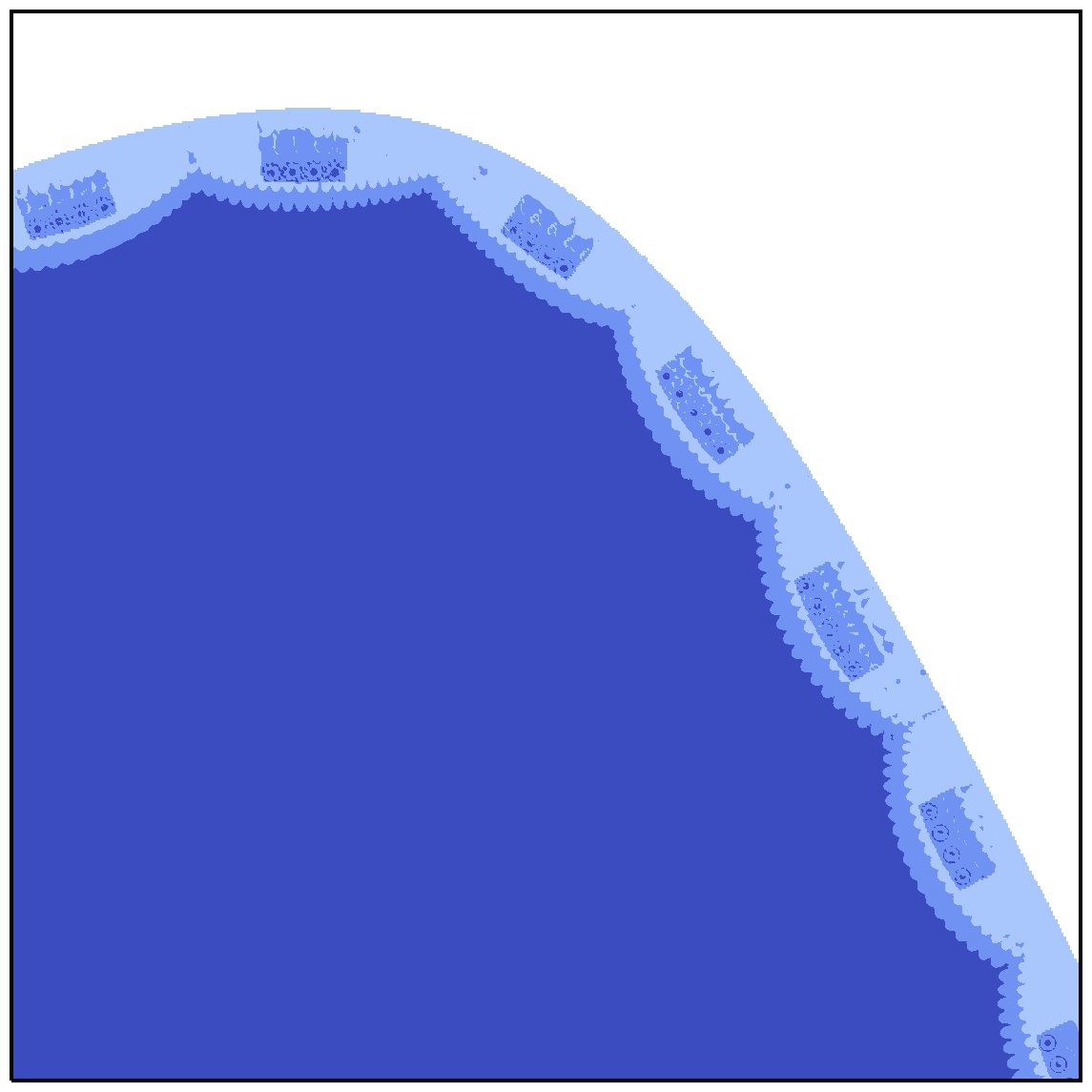}
    \end{minipage}
    \begin{minipage}{0.41\linewidth}
      \includegraphics[width=\textwidth]{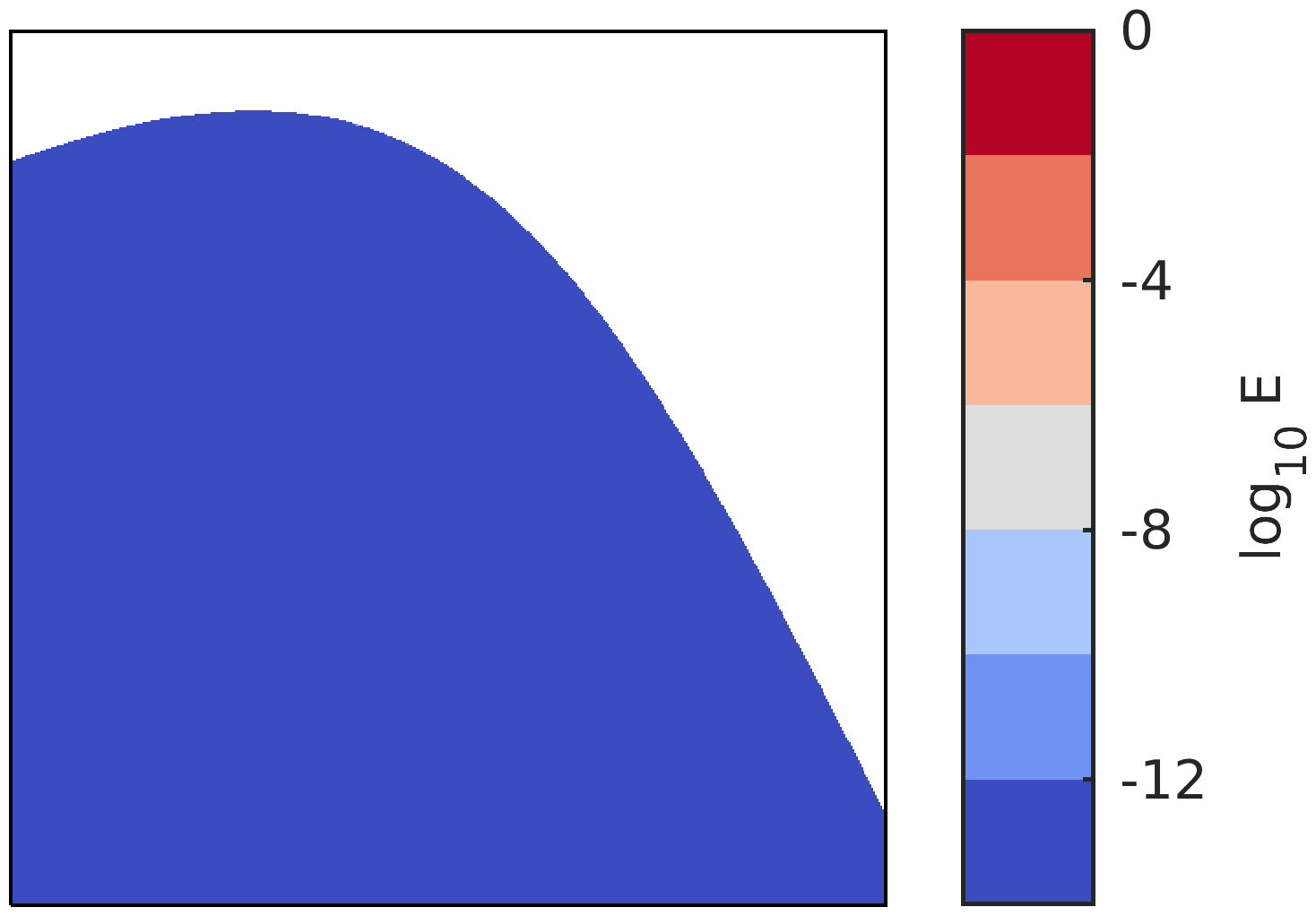}
    \end{minipage}
    \caption{ Errors when evaluating the potential of
      \cref{fig:helmholtz_sol} on a $500 \times 500$ grid, in the region
      highlighted in \cref{fig:helmholtz_err}. The AQBX tolerance is set
      to $10^{-4}$, $10^{-8}$ and $10^{-12}$, from left to right, and
      the region where AQBX is activated is determined using the error
      estimate for the potential \eqref{eq:helmholtz_err}. }
    \label{fig:corrected_field}
  \end{hlbox}
\end{figure}

\begin{figure}[htbp]
  \begin{hlbox}
    \centering
    \begin{minipage}{0.32\textwidth}
      \includegraphics[width=\textwidth]{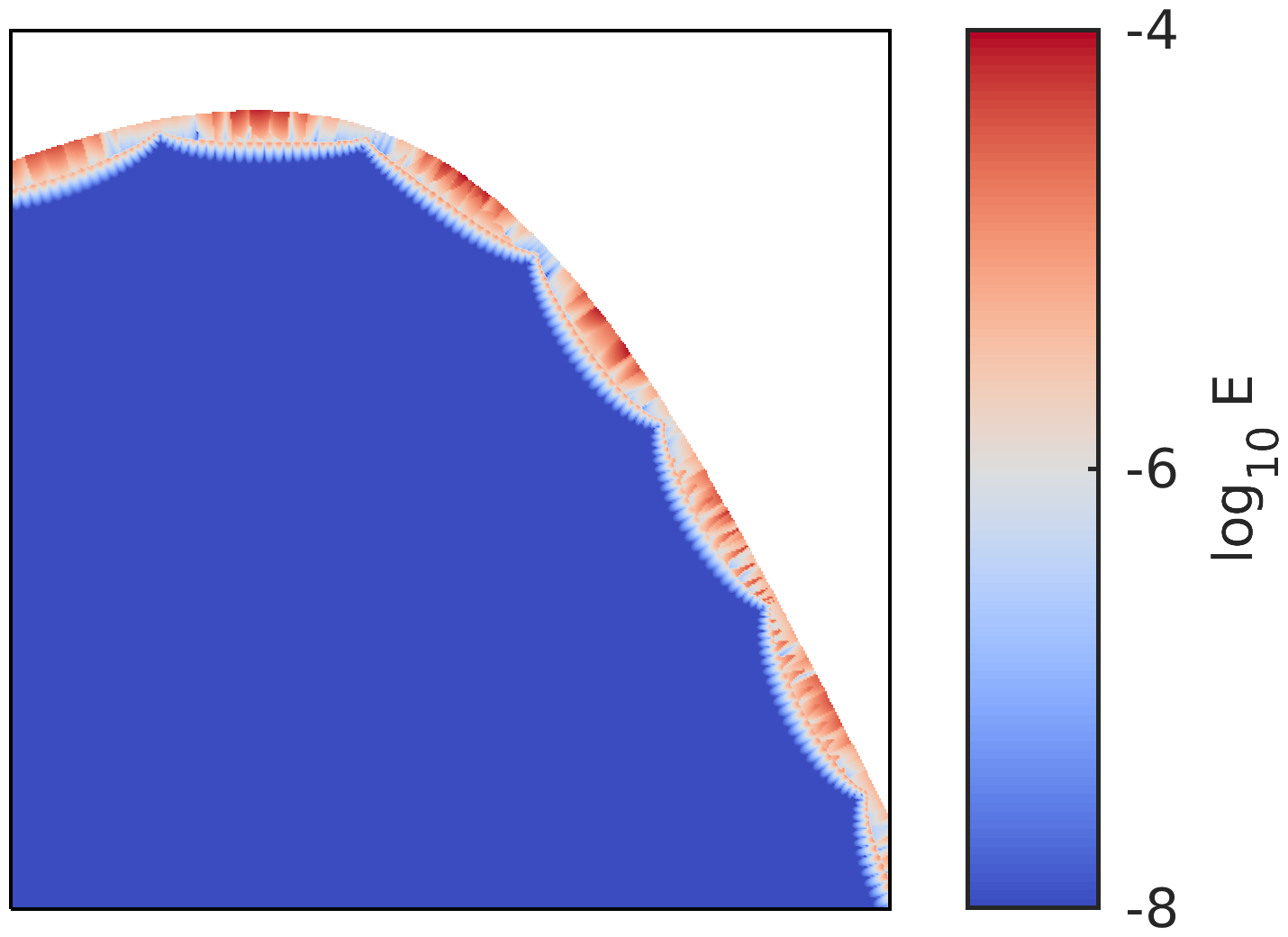}
    \end{minipage}
    \begin{minipage}{0.32\linewidth}
      \includegraphics[width=\textwidth]{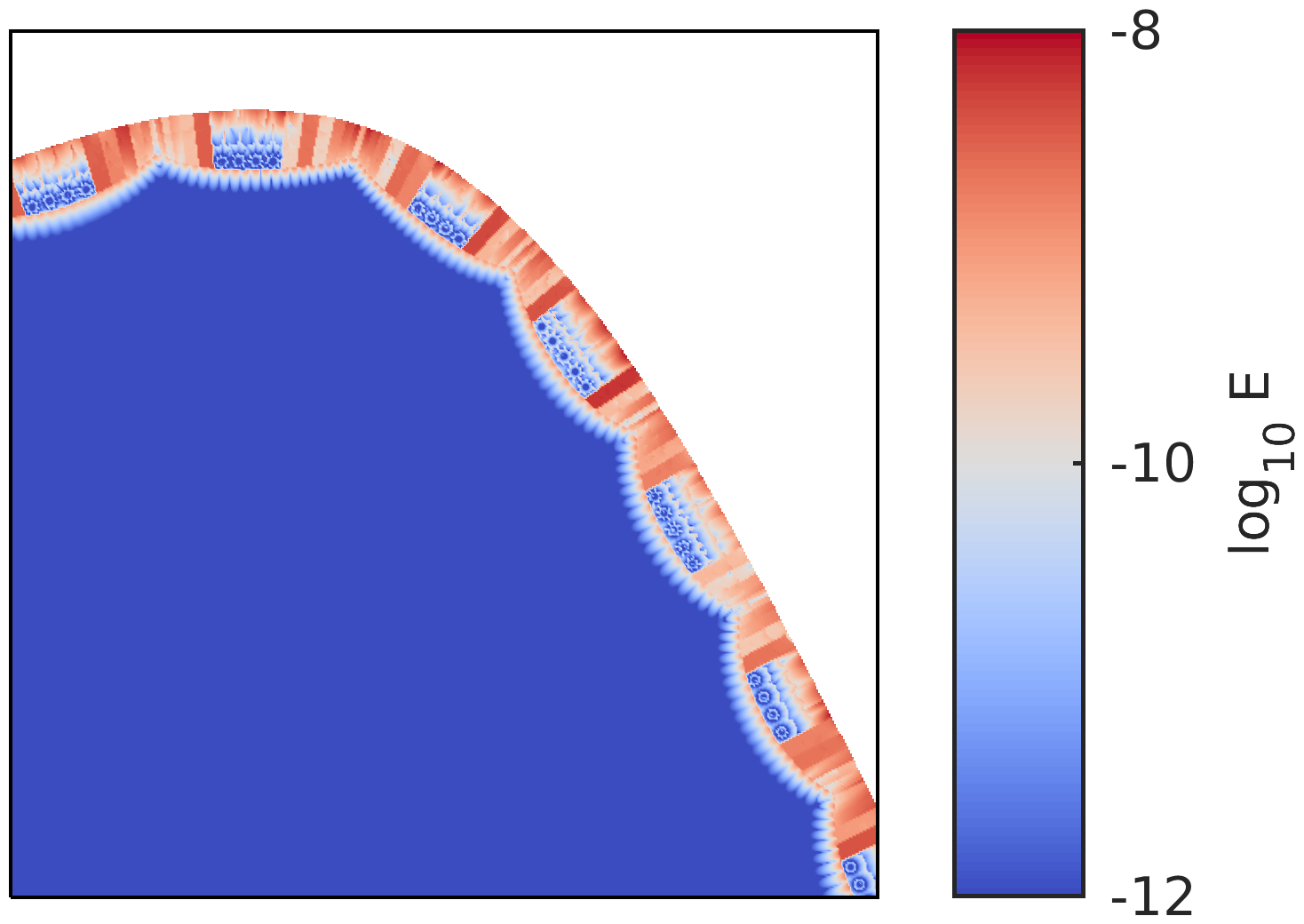}
    \end{minipage}
    \begin{minipage}{0.32\linewidth}
      \includegraphics[width=\textwidth]{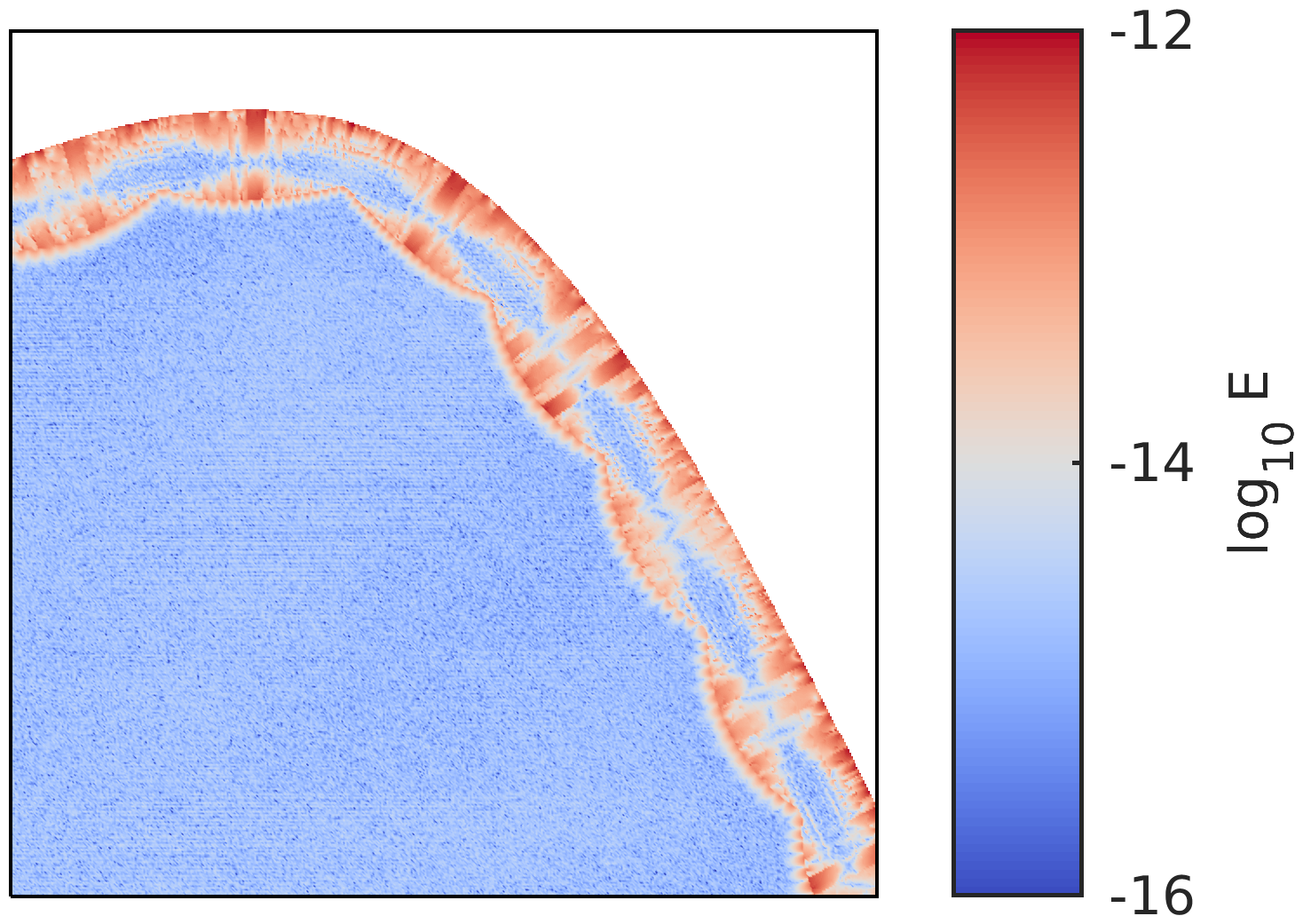}
    \end{minipage}
    \caption{Here the same data is plotted as in
      \cref{fig:corrected_field}, but with compressed colorbars. This
      makes it possible to see the structure of the remaining errors for
      each tolerance.}
    \label{fig:corrected_field_new}
  \end{hlbox}
\end{figure}

\hl{We first perform a few experiments that illustrate how AQBX
  works. \mbox{\Cref{fig:error}} shows the error along the entire
  boundary when evaluating the solution using AQBX and a tolerance of
  $10^{-10}$.} One can clearly see how the magnitude of the
coefficients $a_m$ provides a good overestimate of the truncation
error, while the coefficient error is closely tracked by the estimate
$\estq$. \hl{\mbox{\Cref{fig:corrected_field,fig:corrected_field_new}}}
show example results from when AQBX is used for evaluating the
potential in the domain, where the integral is nearly singular. It can
be seen that AQBX is only activated at the points where it is needed,
and that the potential is then evaluated to the desired accuracy at
those points.

\hl{The lowest error that we can achieve using QBX (both adaptive and
  direct) for this problem is around $10^{-12}$, which can also be
  seen in \mbox{\cref{tab:tol_work}}. This presents a loss of two
  digits of accuracy compared to the error of $10^{-14}$ achieved
  using the direct quadrature away from the boundary. While we can not
  immediately explain this loss, it is reminiscent of the results in
  \mbox{\cite[Tab. 1]{Klockner2013}}. There, the error when using QBX
  for evaluating the double layer is two orders of magnitudes larger
  than for the single layer.}

\subsection{Comparison to direct QBX}
\label{sec:comp-direct-qbx}

We believe that the main benefit of using AQBX rather than a direct
QBX implementation is that the parameter choice is greatly simplified;
given an expansion distance $r$ and a tolerance $\tol$, the upsampling
rate $\upsamp$ and expansion order $p$ are set on the fly as needed at
each expansion center. A second benefit is that setting $\upsamp$ and
$p$ on the fly can save some work, compared to using fixed values
everywhere. In an attempt to quantify this, we now introduce a measure
of the work ($\work$) needed to form a local expansion, in terms of
source evaluations per original source point. If direct QBX is used
with order $p$ and a fixed upsampling rate $\upsamp$, then the work is
given by
\begin{align}
    \work_{\text{QBX}} &= p \upsamp.
\end{align}
If AQBX is used to compute $p$ coefficients, with upsampling rate
$\upsamp_m$ used to evaluate the $m$th coefficient, then the work is
given by
\begin{align}
  \work_{\text{AQBX}} &= \sum_{m=1}^p \upsamp_m .
\end{align}

As a comparison, in \cref{tab:tol_work,tab:rh_work} we measure the
work when computing all expansion coefficients in our reference
problem using AQBX, and compare that to the work needed if $p$ and
$\upsamp$ are fixed everywhere to the minimum values required to
achieve the same accuracy as AQBX. These fixed values are tuned by
hand to the optimal values for this specific problems, \hl{but} our
algorithm \hl{still} gives a slight speedup in our definition of
work. More importantly, our algorithm is in most cases able to keep
the error at the desired order of magnitude without any manual
intervention\hl{. Hand-tuning parameters is on the other hand strictly
  unfeasible in real applications, and will generally result in an
  unnecessarily conservative choice of parameters.}

Our measure of work does not take into account the extra effort needed
to evaluate the estimate of the AQBX scheme. The reported speedup
should therefore be viewed as an upper limit \hl{compared to the
  optimal parameter set}, and as an indicator that automatic parameter
selection does not necessarily have to be more expensive than using
fixed parameters. To minimize the overhead of the scheme one can
evaluate the error estimates \cref{eq:am_err,eq:am_err_helmholtz}
recursively, and the multiple levels of upsampling can for each panel
be computed and reused as needed using a caching algorithm.

\begin{table}[htbp]
  \begin{hlbox}
  \centering
  \begin{tabular}{r|r|r|r|r|r}
    \hline
    Tol. & Eval. error & Exp. order $p$ & Upsamp. rate $\upsamp$ & Work $\work$ &Speedup (avg)\\
    $\tol$ & $\norm{u_p - u}_{\infty}$ & avg (opt) & avg (opt) & avg (opt) & $\work_{\text{QBX}} / \work_{\text{AQBX}}$  \\
    \hline
    $10^{-4~}$ & $1.4\cdot 10^{-4~}$  & ~5.6  (~~4) & 1.1~ (~~2)  & ~6.0 (~~8)  &  1.3 \\ 
    $10^{-6~}$ & $1.7\cdot 10^{-6~}$  & ~7.0  (~~5) & 1.5~ (~~2)  & 10.4 (~10)  &  1.0 \\ 
    $10^{-8~}$ & $1.5\cdot 10^{-8~}$  & ~8.8  (~~7) & 1.9~ (~~3)  & 17.0 (~21)  &  1.2 \\ 
    $10^{-10}$ & $2.2\cdot 10^{-10}$  & 10.3  (~~9) & 2.3~ (~~3)  & 23.2 (~27)  &  1.2 \\ 
    $10^{-12}$ & $2.0\cdot 10^{-12}$  & 12.2  (~11) & 2.6~ (~~4)  & 32.2 (~44)  &  1.4 \\ 
    $10^{-13}$ & $1.1\cdot 10^{-12}$  & 13.3  (~11) & 2.8~ (~~4)  & 37.6 (~44)  &  1.2 \\ 
    \hline
  \end{tabular}
  \caption{Results for varying tolerance and $r/h=1/4$, comparing AQBX
    and direct QBX on our reference problem. Error is measured on
    $\bdry$. Reported parameters are for AQBX an average over all
    expansion centers (avg), and for direct QBX the optimal fixed
    values used at all centers (opt). Note that the smallest error
    that we can obtain is around $10^{-12}$. This limitation holds
    also when using direct QBX.}
  \label{tab:tol_work}
\end{hlbox}
\end{table}

\begin{table}[htbp]
  \begin{hlbox}
  \centering
  \begin{tabular}{r|r|r|r|r|r}
    \hline
    Dist. & Eval. error & Exp. order $p$ & Upsamp. rate $\upsamp$ & Work $\work$ &Speedup (avg)\\
    $r/h$ & $\norm{u_p - u}_{\infty}$ & avg (opt) & avg (opt) & avg (opt) & $\work_{\text{QBX}} / \work_{\text{AQBX}}$  \\
    \hline
    $0.10$ & $1.7\cdot 10^{-10}$  & ~8.0  (~~6) & 4.1~ (~~7)  & 32.8 (~42)  &  1.3 \\ 
    $0.25$ & $2.2\cdot 10^{-10}$  & 10.3  (~~9) & 2.3~ (~~3)  & 23.2 (~27)  &  1.2 \\ 
    $0.50$ & $6.2\cdot 10^{-9~}$  & 13.8  (~12) & 1.7~ (~~2)  & 23.0 (~24)  &  1.0 \\ 
    $0.75$ & $5.4\cdot 10^{-10}$  & 17.7  (~20) & 1.5~ (~~2)  & 27.1 (~40)  &  1.5 \\ 
    $1.00$ & $1.1\cdot 10^{-9~}$  & 22.0  (~30) & 1.8~ (~~2)  & 40.3 (~60)  &  1.5 \\ 
    \hline  
  \end{tabular}
  \caption{Results for varying $r/h$ and tolerance $\tol=10^{-10}$,
    computed in the same way as in \cref {tab:tol_work}.}
  \label{tab:rh_work}
\end{hlbox}
\end{table}

\subsection{Source points close to the boundary}
\label{sec:source-points-close}

\begin{figure}[htbp]
  \begin{hlbox}
  \centering
  \includegraphics[height=0.3\textwidth]{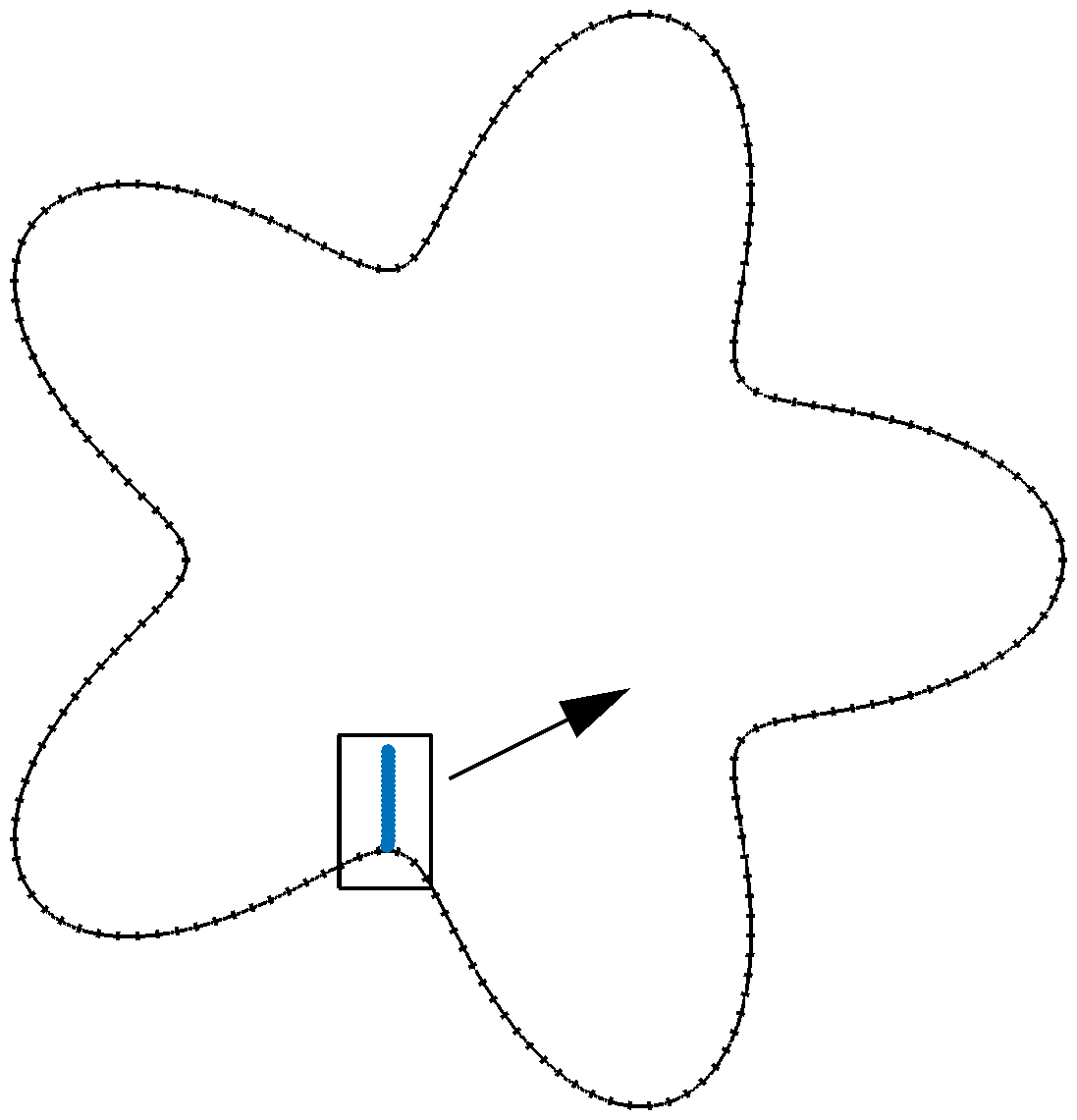}
  \llap{\raisebox{0.05\textwidth}{\includegraphics[height=0.2\textwidth]{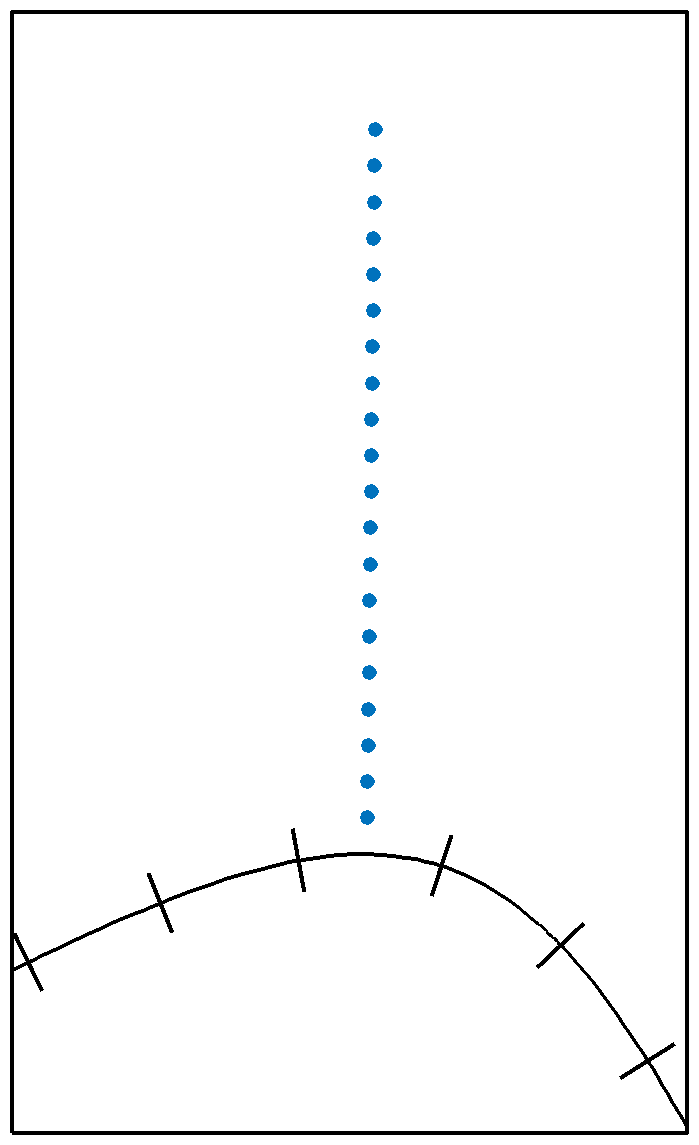}}}
  \includegraphics[width=0.69\textwidth, height=0.3\textwidth]{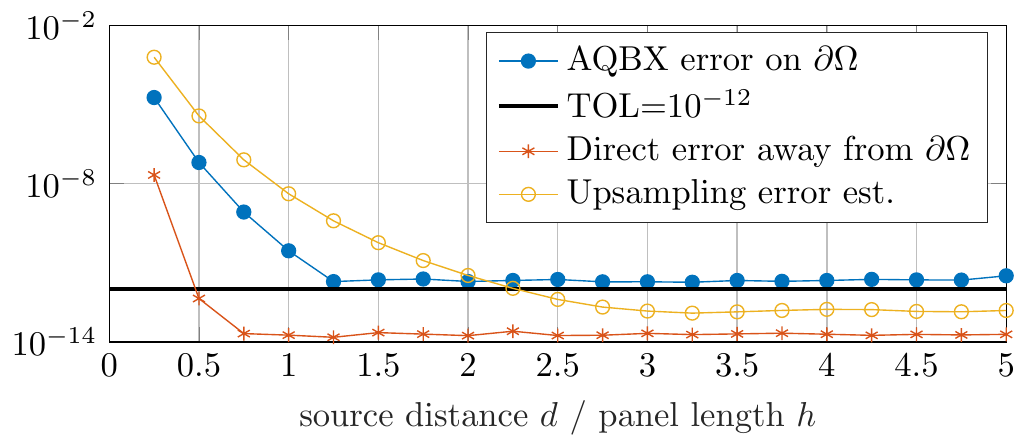}  
  \caption{Results when using the same discretization as in previous cases (200 panels,
    $\wnum=2/h=44.36$), but with Dirichlet boundary conditions given by a point source at
    a distance $d$ from the boundary. We run 20 different cases, with the point source in
    each case given by one of the blue dots in the left graphic. The source strength is
    for each case normalized such that $\norm{u}_\infty=1$ on $\bdry$. AQBX is evaluated
    using a tolerance of $10^{-12}$. The direct error is the maximum relative error in the
    solution when evaluated using direct quadrature on a circle of radius 2 (the outer
    radius of the starfish is 1.3). The upsampling error estimate is the maximum of
    $|\hat\den_{15}| / \norm{\den}_\infty$ over all panels.}
  \label{fig:close_source}
\end{hlbox}
\end{figure}

\hl{In order to push the limits of the algorithm, we run a sequence of
  successively harder test cases. In each case, the Dirichlet boundary
  condition is given by the potential from a source point that lies at
  a distance $d$ from the boundary, with $d$ becoming successively
  smaller. Geometry and results are shown in
  \mbox{\cref{fig:close_source}}. Setting the AQBX tolerance to
  $10^{-12}$, which is the lowest that we can achieve, we see that the
  algorithm starts losing accuracy once the source point is within a
  distance $d=h$ from the boundary. Meanwhile, the error when using
  direct quadrature, measured away from the boundary, does not start
  growing until $d=h/2$.}

\hl{The computation of the QBX coefficients relies on upsampling the
  density $\den$ on each panel using polynomial interpolation. If the
  density is not well represented by its interpolant, then that will
  limit the accuracy of the QBX coefficients.  One way of
  approximating the accuracy of this interpolation is by considering
  the relative magnitude of the highest-order coefficient in the
  Legendre expansion \mbox{\eqref{eq:legendre_interpolant}} of the
  density on each panel, $|\hat\den_{15}| / \norm{\den}_\infty$. We
  take the maximum of this measure over all panels, and denote it the
  ``upsampling error estimate'', shown in
  \mbox{\cref{fig:close_source}}. While not a strictly defined error
  measure, this quantity gives us some indication of how well the
  density is resolved on the grid. It is clear from the figure that
  this error estimate starts growing for source points closer than
  $3 h$ from the boundary. Additionally, it seems that the AQBX error
  grows at the same rate for $d\le h$, though we can not say for
  certain that this is the mechanism governing the AQBX error.}

\hl{The results for this particular test case are encouraging, since
  AQBX does not exhibit any adverse behavior due to the nearby
  singularities. Instead, the error increase appears to be due to a
  lack of resolution.}

\subsection{Solving the integral equation}
\label{sec:solv-integr-equat}

\hl{The above tests indicate that AQBX works well for evaluating layer
  potentials, both close to and on the boundary. Since the method is
  accurate on the boundary, it can also be used to solve the
  underlying integral equation \mbox{\eqref{eq:generic_inteq}}. AQBX
  is then used to evaluate the left-hand side of the discretized
  integral equation (i.e. the matrix-vector product), and a solution
  is found iteratively using GMRES \mbox{\cite{Saad1986}}. The details
  of this are for direct QBX discussed to some extent in
  \mbox{\cite{Klockner2013}}. There, they recommend that the principal
  value integral of the double layer is computed using an average of
  QBX expansions on both sides of the boundary. This was further
  studied in \mbox{\cite{Rahimian2016}}, where they found that
  two-sided expansions had better convergence properties. We here
  follow this recommendation.}

\hl{\mbox{\Cref{tab:inteq_solve}} shows the results when solving our
  test problem for a range of tolerances. We find that the tolerance
  set for the GMRES iterations is matched by both the error in the
  solution $u$ and the smoothness of the density $\den$.  However, we
  find that the AQBX tolerance must be set two orders of magnitude
  smaller than the GMRES tolerance, otherwise GMRES stagnates. We
  believe that this can partly be explained by the fact that AQBX does
  not strictly satisfy its given tolerance. Though beyond the scope of
  the present paper, it would be interesting to further study the
  interplay between the tolerances of AQBX and GMRES. In particular,
  we believe that AQBX could be successfully combined with the Inexact
  GMRES method \mbox{\cite{Simoncini2003}}, which is designed to work
  with a matrix-vector product that has been deliberately made
  inexact.}

\begin{table}[htbp]
  \begin{hlbox}
    \centering
    \begin{tabular}{r|r|r|r|r}
      \hline
      GMRES tol. & AQBX tol. & Iter. count & Rel. error & Upsamp. err. est. \\
      $\tol_\text{GMRES}$      & $\tol_\text{AQBX}$ & AQBX (ref)   & $\frac{\norm{u_p-u}_\infty }{ \norm{u}_\infty}$  & $\max \frac{|\hat\den_{15}| }{ \norm{\den}_\infty}$\\
      \hline
      $10^{-2~}$ & $10^{-4~}$ & ~5 (~~5) & $9.0\cdot 10^{-03}$ & $7.1\cdot 10^{-03}$ \\ 
      $10^{-4~}$ & $10^{-6~}$ & 11 (~11) & $7.1\cdot 10^{-05}$ & $2.7\cdot 10^{-05}$ \\ 
      $10^{-6~}$ & $10^{-8~}$ & 17 (~17) & $5.6\cdot 10^{-07}$ & $3.5\cdot 10^{-07}$ \\ 
      $10^{-8~}$ & $10^{-10}$ & 22 (~22) & $9.0\cdot 10^{-09}$ & $4.6\cdot 10^{-09}$ \\ 
      $10^{-10}$ & $10^{-12}$ & 28 (~28) & $6.4\cdot 10^{-11}$ & $7.8\cdot 10^{-11}$ \\ 
      $10^{-12}$ & $10^{-14}$ & 34 (~33) & $3.9\cdot 10^{-13}$ & $7.9\cdot 10^{-13}$ \\ 
      \hline
    \end{tabular}  
    \caption{Results when solving the problem of \cref{fig:helmholtz_sol} using AQBX, with
      200 panels. The iteration count is the number of GMRES iterations needed to
      converge, and the reference value is the count when using the method of
      \cite{Helsing2015}. The error in the solution is measured on a circle of radius 2,
      and the upsampling error estimate is reported as the maximum over all panels.}
    \label{tab:inteq_solve}
  \end{hlbox}
\end{table}

\FloatBarrier

\section{Conclusions}

We have in this paper formulated a scheme for adaptive quadrature by
expansion (AQBX), which allows for the evaluation of singular and
nearly singular layer potential integrals on a curve discretized using
composite Gauss-Legendre quadrature. The scheme automatically sets
parameters in order to satisfy a given tolerance. This is a
simplification compared to the original QBX scheme, which has a large
parameter space. Given a target tolerance, the only free parameter is
here the expansion radius $r$. Since the remaining parameters are set
on the fly, varying $r$ will mainly affect the speed of the
algorithm. The optimal value for $r$ with respect to speed will be
implementation-dependent, though values of $W$ in \cref{tab:rh_work}
suggest that $r/h$ in the range $0.25$--$0.50$ is a good choice.

The key component of our scheme is the ability to accurately estimate
the magnitude of the quadrature errors in the QBX coefficients. To do
this we have built on the results of \cite{AfKlinteberg2016quad},
where such estimates were reported for a flat panel. We have extended
these to curved panels by taking into account the mapping between a
flat panel and a curved panel. This mapping can be locally constructed
using only the locations of the quadrature nodes on each panel, and
therefore requires no additional analytical information. A side
benefit of our estimates is that the nearly singular quadrature error
of the underlying layer potential can be accurately estimated, which
provides an excellent criterion for when to activate special
quadrature also when other methods are used
\cite{Helsing2008,Barnett2015}. This could also prove useful for QBX
schemes where the expansion is formed by multiple layer potential
evaluations in a neighborhood of the expansion center
\cite{Askham2017,Rahimian2016}.

The focus on a target accuracy is, in our experience, uncommon in the
context of singular and nearly singular quadrature. We do however
believe that this is important if integral equation methods are to be
used in large-scale simulations, where the focus is on achieving a
target accuracy at the lowest possible computational cost.

While several excellent special quadrature methods exist in two
dimensions, methods in three dimensions have not yet reached the same
maturity.  The QBX method has been successfully used on simple
geometries in 3D \cite{AfKlinteberg2016qbx}, while development for
more general use is ongoing. If accurate parameter selection is
important in 2D, it is absolutely essential in 3D, as costs are higher
and the impact of suboptimal parameter choices more severe.  The
principles of the present 2D scheme can be extended to three
dimensions. Estimates for the coefficient errors in 3D were in
\cite{AfKlinteberg2016quad} developed for the special case of
spheroidal surfaces. Developing estimates for general surfaces in 3D
is a topic of ongoing work, and results will be reported at a later
date.

\section*{Acknowledgments}

The authors wish to thank Prof. Johan Helsing for providing an
implementation of the singular integration scheme of
\cite{Helsing2015}.

\bibliography{library}
\bibliographystyle{siamplain_mod}

\end{document}